\documentclass[11pt,a4paper]{article}
\usepackage[T1]{fontenc}
\usepackage[utf8]{inputenc}
\usepackage{authblk}
\usepackage{geometry}
\author[1]{Yue Feng}
\affil[1]{Department of Mathematics, National University of Singapore, Singapore 119076, Singapore}
\author[2]{Ying Ma}
\affil[2]{Beijing Computational Science Research Center, Beijing 100193, China}

\date{}

\usepackage{amssymb,amsfonts,amsmath}
\usepackage{graphicx}
\usepackage{color}
\usepackage{booktabs}
\usepackage{multirow}
\usepackage{threeparttable}
\usepackage{caption}
\usepackage{tabularx}
\usepackage{graphicx}
\usepackage{epsfig}
\usepackage{epstopdf}
\usepackage{subfigure}
\usepackage{multirow}
\usepackage{paralist}
\usepackage{stmaryrd}
\usepackage{booktabs}
\usepackage{amscd}
\usepackage{cite}
\usepackage{amssymb}
\usepackage{verbatim}
\usepackage{multirow}
\usepackage{color}
\usepackage{ntheorem}
\usepackage[all,cmtip]{xy}
\newtheorem{theorem}{Theorem}
\newtheorem{lemma}{Lemma}

\newtheorem{corollary}{Corollary}
\numberwithin{equation}{section}
\numberwithin{theorem}{section}
\numberwithin{lemma}{section}
\numberwithin{remark}{section}
\numberwithin{corollary}{section}

\begin{document}

\title{Error bounds of fourth-order compact finite difference methods for the Dirac equation in the massless and nonrelativistic regime}

\maketitle

\begin{abstract}
We establish the error bounds of fourth-order compact finite difference (4cFD) methods for the Dirac equation in the massless and nonrelativistic regime, which involves a small dimensionless parameter $0 < \varepsilon \le 1$ inversely proportional to the speed of light. In this regime, the solution propagates waves with wavelength  $O(\varepsilon)$ in time and $O(1)$ in space, as well as with the wave speed $O(1/\varepsilon)$ rapid outgoing waves. We adapt the conservative and semi-implicit 4cFD methods to discretize the Dirac equation and rigorously carry out their error bounds depending explicitly on the mesh size $h$, time step $\tau$ and the small parameter $\varepsilon$. Based on the error bounds, the $\varepsilon$-scalability of the 4cFD methods is $h = O(\varepsilon^{1/4})$ and $\tau = O(\varepsilon^{3/2})$, which not only improves the spatial resolution capacity but also has superior accuracy than classical second-order finite difference methods. Furthermore, physical observables including the total density and current density have the same conclusions. Numerical results are provided to validate the error bounds and the dynamics of the Dirac equation with different potentials in 2D is presented.
\end{abstract}


\section{Introduction}
\label{intro}
Since the first proposition by Paul Dirac in 1928, the Dirac equation has been one of the fundamental equations in quantum mechanics \cite{Dirac1,Dirac2}. It describes the motion of elementary spin-1/2 massive particles, such as positrons, electrons, quarks, etc \cite{HPA,NGPNG,SGMH}. As a relativistic wave equation which plays an important role in particle physics, the Dirac equation successfully unifies the principles of quantum mechanics and the theory of special relativity \cite{ACD}. The negative energy solution of the Dirac equation predicted the existence of antimatter, which was previously unobserved and experimentally confirmed by the discovery of the positron a few years later. Recently, due to the development of theoretical studies and experimental explorations for the dynamical properties and structures of graphite and graphene \cite{CG,SGMH,NGMJKGDF}, intense laser-molecule interaction \cite{FLB,FLB1}, as well as 2D and 3D topological insulators \cite{QZ}, the study of the Dirac equation has attracted extensive attention and research interests from numerous scholars.

In this paper, we consider the Dirac equation in the massless and nonrelativistic regime \cite{BY}, which means that the mass of the particle is much less than the mass unit and the wave speed is much less than the speed of light. Similar to the techniques used in \cite{BCJT,BY}, in one dimension (1D) and two dimensions (2D), the Dirac equation in the massless and nonrelativistic regime on the unit torus $\mathbb{T}^d$ ($d = 1, 2$) could be expressed in the two-component form as
\begin{equation}
\label{eq:Dirac_21}
i\partial_t\Phi =  \left(- \frac{i}{\varepsilon}\sum_{j = 1}^{d}
	\sigma_j\partial_j + \frac{1}{\varepsilon}\sigma_3 \right)\Phi+ \left(V(t, \mathbf{x})I_2 - \sum_{j = 1}^{d}A_j(t, \mathbf{x})\sigma_j\right)\Phi, \ {\bf x} \in \mathbb{T}^d,
\end{equation}
where $\Phi : = \Phi(t, {\mathbf{x}}) = (\phi_1(t, {\mathbf{x}}), \phi_2(t, {\mathbf{x}}))^T \in \mathbb{C}^2$ is the complex-valued wave function, $i = \sqrt{-1}$, $\partial_j = \frac{\partial}{\partial x_j}( j = 1, \cdots, d),$ $t$ is time, $\mathbf{x} = (x_1, \cdots , x_d)^T \in \mathbb{T}^d$ is the spatial coordinate vector, $\varepsilon=\frac{v_s}{c}\in (0, 1]$ is a dimensionless parameter inversely proportional to the speed of light, which represents the ratio between the wave velocity and the speed of light. $V:=V(t, \mathbf{x})$ and $\mathbf{A}:=\mathbf{A}(t, \mathbf{x}) = (A_1(t, \mathbf{x}), \cdots, A_d(t, \mathbf{x}))^T $ stand for the real-valued electric potential and magnetic potential, respectively. Besides, $I_2$ is the $2 \times 2$ identity matrix, and $\sigma_1, \sigma_2, \sigma_3$ are the Pauli matrices defined as
	\begin{equation*}
	\label{Pauli}
	\sigma_1 = \begin{pmatrix} 0 &\ \  1\\ 1 &\ \  0 \end{pmatrix}, \quad
	\sigma_2 = \begin{pmatrix} 0 & \ \ -i \\ i &\ \   0\end{pmatrix}, \quad
	\sigma_3 = \begin{pmatrix} 1 &\ \ 0 \\ 0 &\ \  -1 \end{pmatrix}.
	\end{equation*}
In order to study the dynamics of the Dirac equation \eqref{eq:Dirac_21}, the initial data is taken as
\begin{equation}
\label{eq:initial}
\Phi(t=0, \mathbf{x}) = \Phi_0(\mathbf{x}), \quad {\bf x} \in \mathbb{T}^d.
\end{equation}
The Dirac equation \eqref{eq:Dirac_21} is dispersive, time symmetric and conserves the total probability \cite{BCJT}
\begin{equation}
\|\Phi(t, \cdot)\|^2 : = \int_{\mathbb{T}^d} |\Phi(t, \mathbf{x})|^2 d \mathbf{x} = \int_{\mathbb{T}^d} \sum^2_{j=1} |\Phi_j(t, \mathbf{x})|^2 d\mathbf{x} \equiv \|\Phi(0, \cdot)\|^2 = \|\Phi_0\|^2.
\end{equation}
Besides, if the electric and magnetic potentials are both time-independent, i.e., $V(t, \mathbf{x}) = V(\mathbf{x})$ and $A_j(t, \mathbf{x}) = A_j(\mathbf{x}),\,j=1, \cdots, d$, the energy $E(\Phi(t, \cdot)) $ is also conserved
\begin{align}
E(\Phi(t, \cdot)) & : = \int_{\mathbb{T}^d} \Big( -\frac{i}{\varepsilon} \sum_{j=1}^d \Phi^{\ast}\sigma_j\partial_j \Phi + \frac{1}{\varepsilon}\Phi^{\ast}\sigma_3\Phi +V(\mathbf{x}) |\Phi|^2 - \sum_{j=1}^d A_j({\bf x}) \Phi^{\ast}\sigma_j\Phi \Big) d {\bf x} \nonumber\\
	& \equiv E(\Phi(0, \cdot)) = E(\Phi_0), \quad t \geq 0,
\label{eq:Energy}
\end{align}
where $\Phi^{\ast} = {\overline{\Phi}}^T$ with $\overline{\Phi}$ denoting the complex conjugate of $\Phi$.

Introduce the total density $\rho : = \rho(t, \mathbf{x})$ as
\begin{equation}
\label{eq:rho}
	\rho(t, \mathbf{x}) = \sum_{j=1}^2 \rho_j(t, {\bf x}) = \Phi(t, {\bf x})^{\ast}\Phi(t, {\bf x}), \quad {\bf x} \in \mathbb{T}^d,
\end{equation}
with the $j$-th probability density $\rho_j(t, \mathbf{x}) = |\phi_j(t, \mathbf{x}) |^2,\,j=1,2,$ and the current density $\mathbf{J}(t, \mathbf{x}) = (J_1(t, \mathbf{x}),  \cdots, J_d(t, \mathbf{x}))^T$ as
\begin{equation}
\label{eq:J}
J_l(t, \mathbf{x})  = \frac{1}{\varepsilon} \Phi(t, {\bf x})^{\ast} \sigma_l \Phi(t, {\bf x}), \quad l = 1, \cdots, d,
\end{equation}
then the following conservation law could be obtained  \cite{BCJT}
\begin{equation}
\partial_t \rho(t,  \mathbf{x})	+ \nabla \cdot \mathbf{J}(t, \mathbf{x}) = 0, \quad \mathbf{x} \in \mathbb{T}^d, \quad t \geq 0.
\end{equation}

The Dirac equation \eqref{eq:Dirac_21} in different regimes have been widely investigated in the past decades. For the analytical and numerical results in the classical regime, i.e, $\varepsilon = 1$,  we refer to \cite{Das1,Das2,ES,GGT,Gross,SE2,WHJY} and references therein. In the nonrelativistic/semiclassical regime, various numerical methods have been proposed and analyzed including the finite difference time domain (FDTD) methods \cite{BCJT,BHM,MY}, exponential wave integrator Fourier pseudospectral (EWI-FP) method \cite{BCJT,BCJY}, time-splitting Fourier pseudospectral (TSFP) method \cite{BCJY,BCY,BCY11,HJMSZ}, Gaussian bean method \cite{WHJY} and so on \cite{FLB,FLB1,GSX}.

\begin{figure}[ht!]
\begin{minipage}{0.5\textwidth}
\centerline{\includegraphics[width=6.5cm,height=5cm]{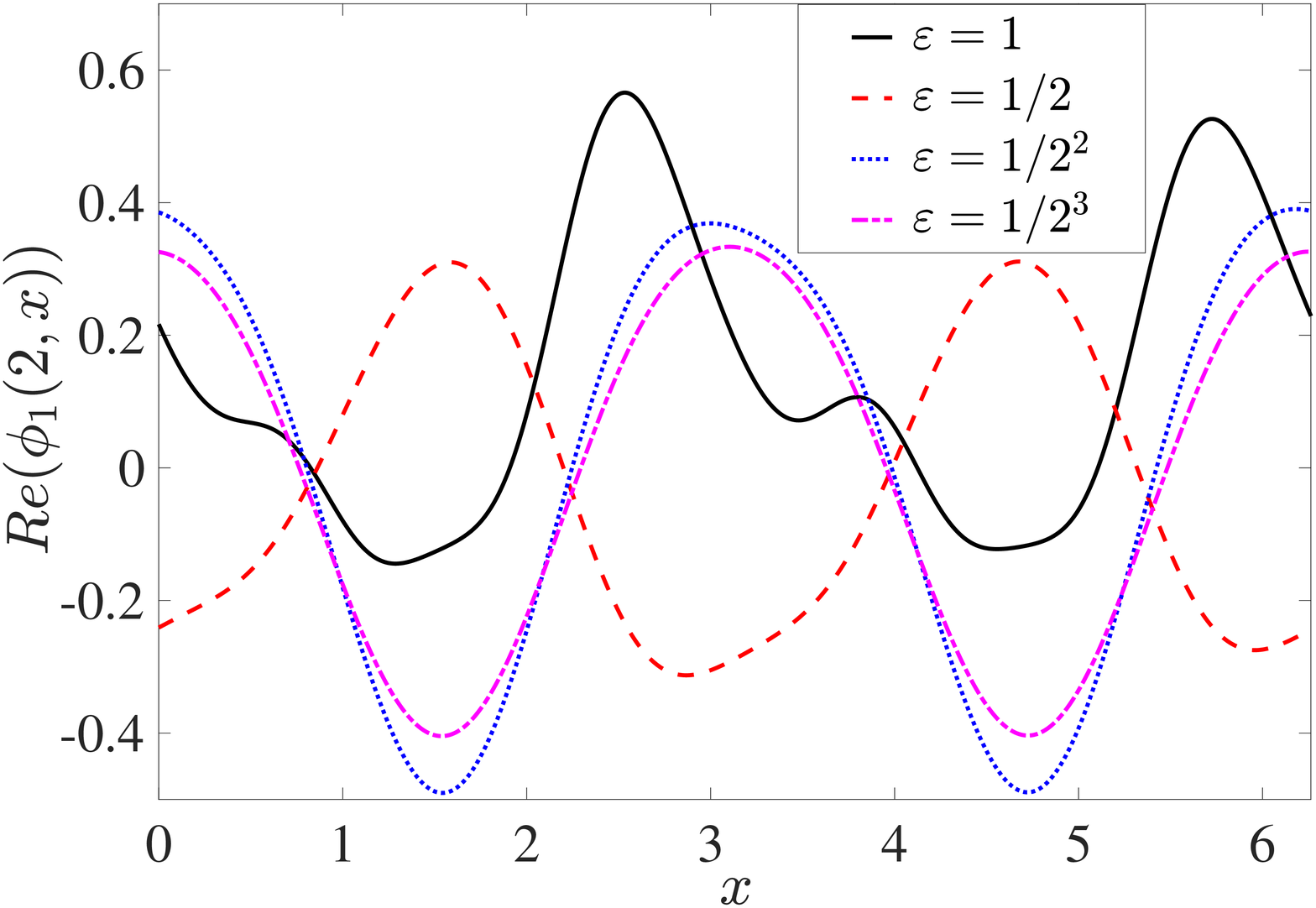}}
\end{minipage}
\begin{minipage}{0.48\textwidth}
\centerline{\includegraphics[width=6.5cm,height=5cm]{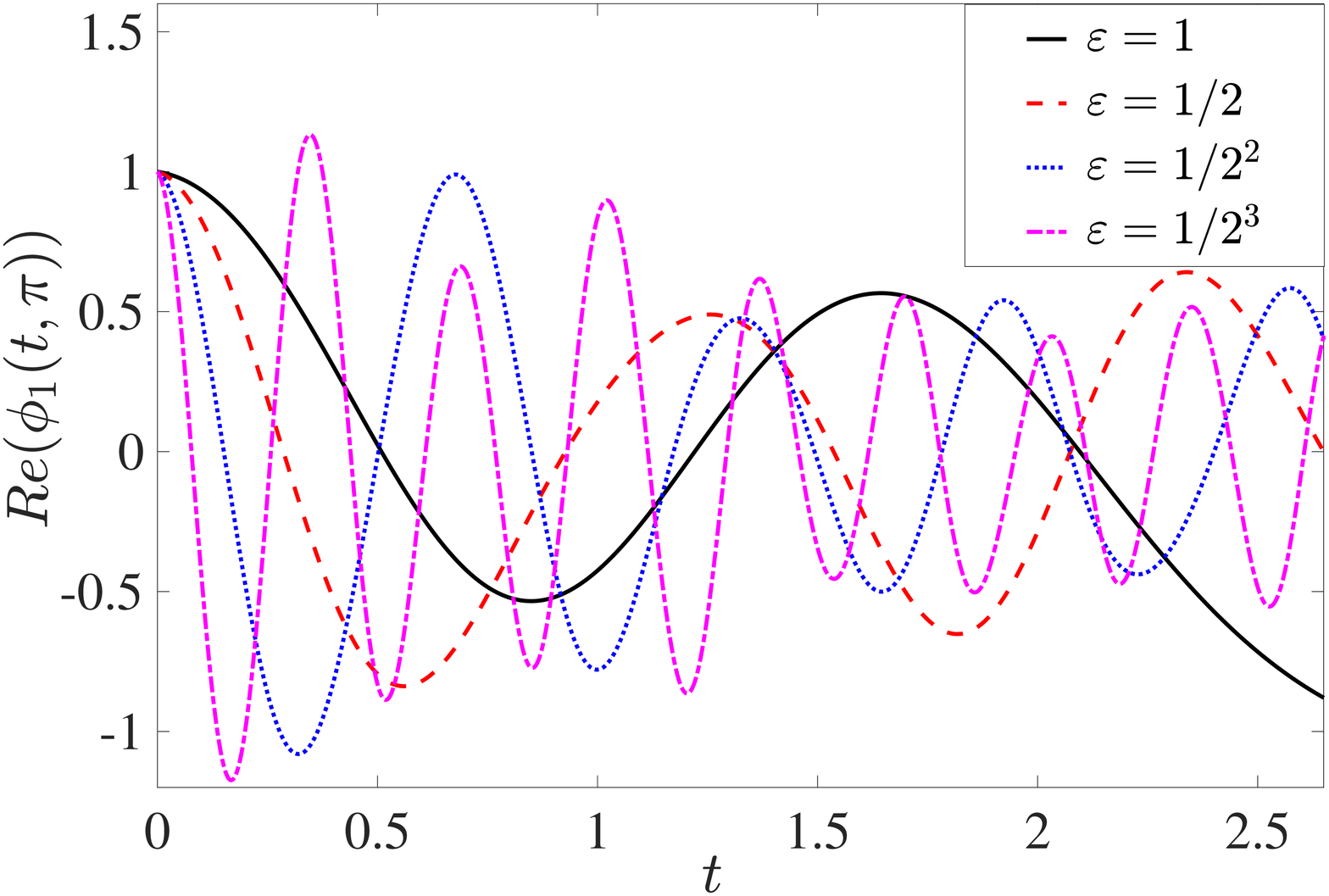}}
\end{minipage}
\caption{The wave function $\phi_1(t=2, x)$ and $\phi_1(t, x=\pi)$ of the Dirac equation \eqref{eq:Dirac_21} with $d=1$ and different $\varepsilon$. Here $\text{Re}(f)$ denotes the real part of $f$}
\label{fig:phi}
\end{figure}

When $0<\varepsilon \ll 1$ in the Dirac equation \eqref{eq:Dirac_21}, i.e., in the massless and nonrelativistic regime, the solution propagates waves with wavelength at $O(1)$ in space and $O(\varepsilon)$ in time. To illustrate this oscillatory nature, Fig. \ref{fig:phi} depicts the solution of the Dirac equation \eqref{eq:Dirac_21} with $d = 1$, $V(t, x) = 1/(1+\sin^2(x))$, $A_1(t, x) = \sin(2x)$ and $\Phi_0(x) = (1/(1+\sin^2(x)), 1/(3+\cos(x)))^T$ for different $\varepsilon$. In fact, in this regime, the highly oscillatory nature of the solution in time brings significant difficulty in numerical simulations and mathematical analysis. Thus, it is very important to adopt effective numerical methods to study the dynamics of the Dirac equation \eqref{eq:Dirac_21} when $\varepsilon \to  0^{+}$ and carry out rigorous error bounds, especially the explicit dependence on the mesh size $h$ and time step $\tau$ as well as the small parameter $\varepsilon$.

High order compact finite difference methods could achieve expected accuracy with less grid points, which are able to improve the spatial resolution capacity especially for $0 <\varepsilon \ll 1$. The fourth-order compact finite difference (4cFD) method is a simple scheme to attain  higher spatial order with the same number of grids for the central difference method \cite{Liao,MAM,ZW}. Recently, the 4cFD method has been used to solve the (nonlinear) Schr\"odinger equation \cite{HJKW,WZ}, Klein--Gorden equation \cite{Feng,LLG}, Dirac equation \cite{LL}, Burgers' equation \cite{Liao} and so on. For more details, we refer to \cite{LCP,LSS,MAM,WMLG,WGX} and references therein.

To the best of our knowledge, there is no study on the the fourth-order compact finite difference (4cFD) schemes for the Dirac equation \eqref{eq:Dirac_21} in the massless and nonrelativistic regime. The aim of this paper is to combine the 4cFD discretization in space with the implicit/semi-implicit temporal discretization to numerically solve the Dirac equation \eqref{eq:Dirac_21} and carry out the error bounds in the massless and nonrelativistic regime. Compared with the second-order finite difference methods in \cite{MY2}, the 4cFD methods not only have higher accuracy in terms of the mesh size $h$ for the fixed $\varepsilon$, but also have better spatial resolution in terms of $\varepsilon$. Based on our rigorous error estimates, in order to get `correct' numerical approximations of the Dirac equation \eqref{eq:Dirac_21}, the $\varepsilon$-scalability of the 4cFD methods should be taken as
\begin{equation*}
h = O(\varepsilon^{1/4}) \quad \mbox{and} \quad \tau = O(\varepsilon^{3/2}), \quad 0 < \varepsilon \leq 1,
\end{equation*}
which performs much better than the classical finite difference schemes for the spatial resolution.

The rest of this paper is organized as follows. In Section \ref{sec:2}, the implicit and semi-implicit 4cFD methods are presented for the Dirac equation \eqref{eq:Dirac_21} with stability conditions analyzed. In Section \ref{sec:3}, error bounds of these two 4cFD schemes for the Dirac equation \eqref{eq:Dirac_21} are rigorously carried out. Numerical results are shown in Section \ref{sec:4} to confirm our error estimates and study the dynamics of the Dirac equation \eqref{eq:Dirac_21} in 2D. Finally, some conclusions are drawn in Section \ref{sec:5}. Throughout this paper, the notation $A \lesssim B$  is used with the meaning that there exists a generic constant $C > 0$ independent of the mesh size $h$ and time step $\tau$ as well as the parameter $\varepsilon$, such that $|A| \leq C B$.

\section{The 4cFD methods and their analysis}
\label{sec:2}
In the section, we adapt the implicit and semi-implicit fourth-order compact finite difference (4cFD) methods to solve the Dirac equation \eqref{eq:Dirac_21} and analyze their mass and energy conservation as well as the stability conditions. For simplicity of presentation, here we only carry out the numerical schemes and corresponding analysis in 1D. It is straightforward to generalize to 2D and the results remain valid without modifications. The Dirac equation \eqref{eq:Dirac_21} in 1D on the computational domain $\Omega = (a, b)$ with periodic boundary conditions collapses to
\begin{align}
\label{eq:Dirac_1D}
&i\partial_t\Phi =  \left(- \frac{i}{\varepsilon} \sigma_1 \partial_x + \frac{1}{\varepsilon}\sigma_3\right)\Phi+ \left(V(t, x)I_2 - A_1(t, x)\sigma_1\right)\Phi, \ x \in \Omega,\ t > 0,\\
\label{eq:ib}
&\Phi(t, a) = \Phi(t, b),\ t \geq 0; \ \Phi(0, x) = \Phi_0(x),\ x \in \overline{\Omega},
\end{align}
where $\Phi := \Phi(t, x)$ and  $\Phi_0(a) = \Phi_0(b)$.

\subsection{The 4cFD methods}
Choose the mesh size $h : = \Delta x = (b - a)/N$ with $N$ being an even positive integer, time step size $\tau := \Delta t > 0$, and denote the grid points and time steps as
\begin{equation}
 x_j :=a +j h,\quad j = 0, 1, \cdots, N; \quad t_n : = n\tau,\quad n = 0, 1, 2, \cdots.
\end{equation}
Denote $X_N = \{U=(U_0, U_1,\cdots,U_N)^T \ | \ U_j \in \mathbb{C}^2, j = 0,1, \cdots, N,  U_0 = U_N\}$ and let $U_{-1} = U_{N-1}$, $U_{N+1} = U_1$ if they are involved. Define the index set $
\mathcal{T}_N = \{l | l = -N/2,-N/2+1,\cdots,N/2-1\}$, and $\mu_l = \frac{2\pi l}{b-a}\in \mathbb{C}^2$ for $l \in \mathcal{T}_N$,
then for any $U \in X_N$, its corresponding Fourier representation is \cite{ST}
\begin{equation}
U_j = \sum_{l \in \mathcal{T}_N}\widetilde{U}_l e^{i\mu_l(x_j-a)} = \sum_{l \in \mathcal{T}_N}\widetilde{U}_l e^{2ijl\pi/N}, \quad j = 0,1, \cdots, N,
\end{equation}
where
\begin{equation}
\widetilde{U}_l = \frac{1}{N} \sum_{j=0}^{N-1} U_j e^{-2ijl\pi/N}, \quad
\widetilde{U}_l \in \mathbb{C}^2,\,l \in \mathcal{T}_N.
\end{equation}
In the space $X_N,$ the $l^2$-norm and $l^{\infty}$-norm are given  as
\begin{equation}
 \|U\|^2_{l^2} = h \sum_{j=0}^{N-1} |U_j|^2,\quad \|U\|_{l^{\infty}} = \max_{0 \leq j\leq N-1} |U_j| ,\quad U \in X_N.
\end{equation}

Let $\Phi^n_j$ be the numerical approximation of $\Phi(t_n, x_j)$, $V^n_j = V(t_n, x_j)$, $V^{n+1/2}_j = V(t_n + \tau/2, x_j)$,  $A^n_{1, j} = A_1(t_n, x_j)$ and $A^{n+1/2}_{1, j} = A_1(t_n+\tau/2, x_j)$ for $0 \leq j \leq N$ and $n \geq 0$. Denote $\Phi^n = (\Phi^n_0, \Phi^n_1, \cdots, \Phi^n_N)^T \in X_N$ as the solution vector at $t = t_n$. Introduce the following finite difference discretization operators
\begin{equation*}
\begin{split}
\delta^+_t \Phi^n_j = \frac{\Phi^{n+1}_j - \Phi^n_j}{\tau},& \quad  \delta_t \Phi^n_j = \frac{\Phi^{n+1}_j - \Phi^{n-1}_j}{2\tau},\quad \delta_x \Phi^n_j = \frac{\Phi^n_{j+1} - \Phi^n_{j-1}}{2h},\\
\mathcal{A}_h\Phi^n_j = & \left(I + \frac{h^2}{6} \delta^2_x\right)\Phi^n_j = \frac{\Phi^n_{j-1}+4\Phi^n_j + \Phi^n_{j+1}}{6},
\end{split}
\end{equation*}
and the average vector
\begin{equation*}
\Phi^{n+\frac{1}{2}}_j	 = \frac{\Phi^{n+1}_j + \Phi^n_j}{2}.
\end{equation*}
It is easy to check that $\mathcal{A}_h(\partial_x \Phi(x_j)) = (\Phi(x_{j+1})-\Phi(x_{j-1}))/{2h} + O(h^4)$. A fourth-order compact approximation is implemented by replacing $\partial_x \Phi(t_n, x_j)$ by $\mathcal{A}^{-1}_h\left(\delta_x \Phi^n_j\right)$ \cite{LL,LSS}.

Combining the fourth-order compact finite difference discretization in space with the implicit/semi-implicit temporal discretization, we have the following two 4cFD schemes for $j = 0, 1, \cdots, N-1$:

I. The implicit 4cFD method
\begin{equation}\label{eq:im4cFD}
i \delta_{t}^{+} \Phi_{j}^{n}=\frac{1}{\varepsilon}\Big(-i \sigma_{1} \mathcal{A}^{-1}_h\delta_{x} +\sigma_{3} \Big) \Phi_{j}^{n+\frac{1}{2}}+\left(V_{j}^{n+\frac{1}{2}} I_{2}-A_{1, j}^{n+\frac{1}{2}} \sigma_{1}\right) \Phi_{j}^{n+\frac{1}{2}},\ n \geq 0.
\end{equation}

II. The semi-implicit 4cFD method
\begin{equation}
\label{eq:4cFD}
i\delta_t \Phi^n_j = \left(-\frac{i}{\varepsilon}\sigma_1\mathcal{A}^{-1}_h\delta_x+ \frac{1}{\varepsilon}\sigma_3\right)\frac{\Phi^{n+1}_j+ \Phi^{n-1}_j  }{2} + \left(V^n_j I_2 -A^n_{1, j}\sigma_1 \right)\Phi^{n}_j,	\ n \geq 1.
\end{equation}
The boundary and initial conditions \eqref{eq:ib} for the 4cFD methods are discretized as
\begin{equation}
\label{eq:bic}
\Phi^{n+1}_0 = \Phi^{n+1}_N,\quad	 \Phi^{n+1}_{-1} = \Phi^{n+1}_{N-1}, \quad n \geq 0;\quad \Phi^0_j = \Phi_0(x_j), \quad j = 0,1,\cdots, N.
\end{equation}
According to the Taylor expansion and the Dirac equation \eqref{eq:Dirac_21}, the first time step for the semi-implicit 4cFD \eqref{eq:4cFD} could be designed as ($j = 0, 1, \cdots, N$)
\begin{equation}
\label{eq:phi1}
\Phi^1_j = \Phi_0(x_j) - \sin\left(\frac{\tau}{\varepsilon}\right)\sigma_1\Phi'_0(x_j) - i\left[\sin\left(\frac{\tau}{\varepsilon}\right)\sigma_3 + \tau\left(V^0_j I_2- A^0_{1, j}\sigma_1\right)\right]\Phi_0(x_j).	\end{equation}
Here, we adopt $\frac{1}{\tau} \sin\left(\frac{\tau}{\varepsilon}\right)$ instead of $\frac{1}{\varepsilon}$ such that \eqref{eq:phi1} is second order in terms of $\tau$ for any fixed $\varepsilon \in (0,1]$ and $\|\Phi^1\|_{l^{\infty}} \lesssim 1$ \cite{BCJT,MY}.

The 4cFD methods are time symmetric, i.e. they are unchanged under $\tau \leftrightarrow -\tau$ and $n+1 \leftrightarrow n$ for the implicit 4cFD method or $n+1 \leftrightarrow n-1$ for the semi-implicit 4cFD method, and their memory cost are both $O(N)$.  The implicit 4cFD method is an implicit scheme in the sense that at each time step for $n \geq 0$, the corresponding linear system is coupled and needs to be solved by an iterative solver or a direct solver. As a consequence, the computational cost per time step is usually much larger than $O(N)$, especially in 2D. The semi-implicit 4cFD method is also an implicit scheme, but at each time step for $n \geq 1$, the corresponding linear system could be decoupled and solved explicitly in the phase (Fourier) space as
\begin{equation}
\widetilde{(\Phi^{n+1})}_l = \left(i \gamma_l I_2 - \frac{\tau\sin(\mu_l h)}{\varepsilon h}\sigma_1-\frac{\tau}{\varepsilon} \gamma_l \sigma_3\right)^{-1} L^n_l,\quad l \in \mathcal{T}_N,	
\end{equation}
where
\begin{align}
\label{eq:gamma}
&\gamma_l = \frac{\cos(\mu_l h)+2}{3}, \\
&L^n_l = \left(i \gamma_l I_2 +  \frac{\tau\sin(\mu_l h)}{\varepsilon h}\sigma_1 + \frac{\tau}{\varepsilon}\gamma_l \sigma_3\right) \widetilde{(\Phi^{n-1})}_l + 2\tau \gamma_l \widetilde{(G^n\Phi^n)}_l,
\end{align}
and $G^n = (G^n_0, G^n_1, \cdots, G^n_N)^T \in X_N$ with $G^n_j = V^n_j I_2 -A^n_{1, j}\sigma_1 $ for $j = 0,1,\cdots, N$. Therefore, its computational cost per time step is $O(N \ln N)$.

\subsection{Mass and energy conservation}
For the implicit 4cFD scheme, we could derive the following conservative properties.
\begin{lemma}
\label{lemma:conserv1}
The implicit 4cFD \eqref{eq:im4cFD} conserves the mass in the discretized level as
\begin{equation}
\label{eq:mass}
\|\Phi^n\|^2_{l^2} := h \sum^{N-1}_{j=0} |\Phi^n_j|^2 \equiv h \sum^{N-1}_{j=0} |\Phi^0_j|^2 = \|\Phi^0\|^2_{l^2},\quad n\ge0.
\end{equation}
Moreover, if the electric potential $V(t, x) = V(x)$ and magnetic potential $A_1(t, x) = A_1(x)$ are both time-independent, the implicit 4cFD \eqref{eq:im4cFD} also conserves the discretized energy as
\begin{align}
\label{eq:energy_con}
E^n_h := \ & h\sum^{N-1}_{j=0} \left[-\frac{i}{\varepsilon} (\Phi^n_j)^{\ast}\sigma_1\mathcal{A}^{-1}_h \delta_x \Phi^n_j + \frac{1}{\varepsilon}(\Phi^n_j)^{\ast}\sigma_3 \Phi^n_j + V_j|\Phi^n_j|^2 - A_{1, j}(\Phi^n_j)^{\ast}\sigma_1 \Phi^n_j \right] \nonumber\\
\equiv \ & E^0_h, \quad n \geq 0,
\end{align}
where $V_j = V(x_j)$ and $A_{1, j} = A_1(x_j)$ for $j = 0,1, \cdots, N$.
\end{lemma}
\emph{Proof.}
For the mass conservation \eqref{eq:mass}, multiplying both sides of the 4cFD scheme \eqref{eq:im4cFD} from the  left by $h\tau (\Phi^{n+1/2}_j)^{\ast}$ and taking the imaginary part, we obtain for $j = 0, 1, \cdots, N-1$,
\begin{equation*}
h|\Phi^{n+1}_j|^2 = h|\Phi^n_j|^2 - \frac{\tau h}{\varepsilon}\left[(\Phi^{n+1/2}_j)^{\ast}\sigma_1 \mathcal{A}^{-1}_h\delta_x \Phi^{n+1/2}_j + (\Phi^{n+1/2}_j)^{T} \sigma_1\mathcal{A}^{-1}_h \delta_x \overline\Phi^{n+1/2}_j\right].	
\end{equation*}
Since $\mathcal{A}^{-1}_h$ is a self-adjoint operator, summing above equations for $j = 0, 1, \cdots, N-1$, we could get
 \begin{align*}
 \|\Phi^{n+1}\|^2_{l^2} = &\ \|\Phi^{n}\|^2_{l^2} - \frac{\tau h}{\varepsilon}\sum^{N-1}_{j=0}\left[(\Phi^{n+1/2}_j)^{\ast}\sigma_1\mathcal{A}^{-1}_h\delta_x \Phi^{n+1/2}_j \right. \\
 &\ + \left. (\Phi^{n+1/2}_j)^{T} \sigma_1\mathcal{A}^{-1}_h \delta_x \overline\Phi^{n+1/2}_j\right]\\
    = & \  \|\Phi^{n}\|^2_{l^2} - \frac{\tau}{2\varepsilon}\sum^{N-1}_{j=0}\left[(\Phi^{n+1/2}_j)^{\ast}\sigma_1\mathcal{A}^{-1}_h  \Phi^{n+1/2}_{j+1} + (\Phi^{n+1/2}_j)^{T} \sigma_1 \mathcal{A}^{-1}_h \overline\Phi^{n+1/2}_{j+1} \right.\\
      &\left. \ - (\Phi^{n+1/2}_{j+1})^{\ast}\sigma_1\mathcal{A}^{-1}_h  \Phi^{n+1/2}_j - (\Phi^{n+1/2}_{j+1})^{T} \sigma_1\mathcal{A}^{-1}_h  \overline\Phi^{n+1/2}_j \right]  \nonumber\\
 = & \ \|\Phi^{n}\|^2_{l^2}, \quad n \geq 0,
\end{align*}
 which implies the mass conservation \eqref{eq:mass} by induction.

For the energy conservation \eqref{eq:energy_con}, multiplying both sides of the 4cFD scheme \eqref{eq:im4cFD} from the left by $2h (\Phi^{n+1}_j - \Phi^n_j)^{\ast}$ and taking the real part, we arrive at
 \begin{align}
 & -\frac{h}{\varepsilon} \text{Re}\left[i(\Phi^{n+1}_j - \Phi^n_j)^{\ast} \sigma_1 \mathcal{A}^{-1}_h \delta_x(\Phi^{n+1}_j+\Phi^n_j)\right] + \frac{h}{\varepsilon}\left[(\Phi^{n+1}_j)^{\ast}\sigma_3 \Phi^{n+1}_j - (\Phi^{n}_j)^{\ast}\sigma_3 \Phi^{n}_j\right] \nonumber \\
 & \quad+  h V_j (|\Phi^{n+1}_j|^2-|\Phi^{n}_j|^2) - h A_{1, j} \left[(\Phi^{n+1}_j)^{\ast}\sigma_1 \Phi^{n+1}_j - (\Phi^{n}_j)^{\ast}\sigma_1 \Phi^{n}_j\right] = 0.
 \end{align}
Similarly, summing above equations for $ j = 0, 1, \cdots, N-1$, and noticing the summation by parts formula and the self-adjoint operator $\mathcal{A}^{-1}_h$, we have
\begin{equation*}
\begin{split}
& \frac{h}{\varepsilon} \sum^{N-1}_{j=0} \text{Re} \left(i(\Phi^{n+1}_j - \Phi^n_j)^{\ast}\sigma_1 \mathcal{A}^{-1}_h\delta_x (\Phi^{n+1}_j +\Phi^n_j)\right) \\
= & \  i \frac{h}{\varepsilon}\sum^{N-1}_{j=0}(\Phi^{n+1}_j)^{\ast}\sigma_1\mathcal{A}^{-1}_h\delta_x\phi^{n+1}_j - i\frac{h}{\varepsilon}\sum^{N-1}_{j=0}(\Phi^{n}_j)^{\ast}\sigma_1\mathcal{A}^{-1}_h\delta_x\Phi^{n}_j,	
\end{split}
\end{equation*}
and
\begin{equation*}
\begin{split}
& \ -i \frac{h}{\varepsilon} \sum^{N-1}_{j=0}\left((\Phi^{n+1}_j)^{\ast}\sigma_1\mathcal{A}^{-1}_h \delta_x \Phi^{n+1}_j -(\Phi^{n}_j)^{\ast}\sigma_1\mathcal{A}^{-1}_h \delta_x\Phi^{n}_j \right)\\
& \ + \frac{h}{\varepsilon}\sum^{N-1}_{j=0}\left((\Phi^{n+1}_j)^{\ast}\sigma_3 \Phi^{n+1}_j - (\Phi^n_j)^{\ast}\sigma_3\Phi^n_j\right) + h \sum^{N-1}_{j=0} V_j (|\Phi^{n+1}_j|^2 - |\Phi^n_j|^2)  \\
& \ -  h \sum^{N-1}_{j=0} A_{1, j}\left((\Phi^{n+1}_j)^{\ast} \sigma_1\Phi^{n+1}_j - (\Phi^n_j)^{\ast}\sigma_1\Phi^n_j\right) = 0,
\end{split}
\end{equation*}
which completes the proof.
\hfill $\square$ \bigskip

\subsection{Stability analysis}
Denote $\Omega_T = [0,T_0] \times \Omega$ and take $V_{\max}:= \max \limits_{(t, x) \in \overline{\Omega}_{T}}|V(t, x)|,$ $ A_{1, \max}:= \max \limits_{(t, x) \in \overline{\Omega}_{T}}|A_1(t, x)|$.
Following the standard von Neumann method \cite{Smith}, we could derive the following results for the 4cFD methods.
\begin{lemma}
(i) For any $0<\varepsilon\le1$ and $h,\,\tau>0$, the implicit 4cFD \eqref{eq:im4cFD} with \eqref{eq:bic} is unconditionally stable.

\noindent
(ii) For any $0 < \varepsilon \leq 1$ and $h > 0$, the semi-implicit 4cFD \eqref{eq:4cFD} with \eqref{eq:bic} and \eqref{eq:phi1} is stable under the stability condition
\begin{equation}
\label{eq:stability}
0 < \tau \leq \frac{1}{V_{\max}+ A_{1, \max}}.
\end{equation}
\end{lemma}
\emph{Proof.}
(i) Plugging
\begin{equation}\label{prlefd1}
\Phi_{j}^{n}=\sum_{l\in\mathcal{T}_N} \eta_{l}^{n} \widetilde{\left(\Phi^{0}\right)}_{l} e^{i \mu_{l}\left(x_{j}-a\right)}=\sum_{l\in\mathcal{T}_N} \eta_{l}^{n} \widetilde{\left(\Phi^{0}\right)}_{l} e^{2 i j l \pi / N},\quad 0 \leq j \leq N,
\end{equation}
with $\eta_l^n\in \mathbb{C}$ being the amplification factor of the $l$-th mode in the phase space and $\widetilde{\left(\Phi^{0}\right)}_{l}$ being the Fourier coefficient at $n=0$, respectively, into the implicit 4cFD \eqref{eq:im4cFD}, and using the orthogonality of the Fourier series, we get for $l\in\mathcal{T_N}$,
\begin{equation}\label{prlefdX4}
\left|\left(\eta_{l}-1\right) I_{2}+
\frac{i \tau}{2} (\eta_{l}+1)\left(\frac{\sigma_{3}}{\varepsilon} +V^{0} I_{2}-A_{1}^{0} \sigma_{1}+\frac{\sin \left(\mu_{l} h\right)}{\varepsilon h\gamma_l} \sigma_{1}\right)\right|=0,
\end{equation}
with $\gamma_l$ defined in \eqref{eq:gamma}. Denoting
\begin{equation}
 \theta_{l}=-V^{0} \pm \frac{1}{\varepsilon h} \sqrt{h^{2}+\left(-\varepsilon A_{1}^{0} h+ \frac{\sin \left(\mu_{l} h\right)}{\gamma_l}\right)^{2}},	
\end{equation}
we could obtain
\begin{equation}\label{prcnfdX5}
\eta_{l} = \frac{2+i \tau \theta_{l}}{2-i \tau \theta_{l}}, \quad l \in \mathcal{T}_N,
\end{equation}
which immediately implies $\left|\eta_{l}\right| =1 $ for $l \in \mathcal{T}_N$, so  the implicit 4cFD method is unconditionally stable.

(ii) Since the implicit part $\left(-\frac{i}{\varepsilon}\sigma_1\mathcal{A}_h^{-1}\delta_x+
\frac{1}{\varepsilon}\sigma_{3}\right) \frac{\Phi_{j}^{n+1}+\Phi_{j}^{n-1}}{2}$ is automatically stable, here we just need to focus on the explicit part $i\delta_t\Phi_j^n=\big(V_{\max}I_2-A_{1,\max}\sigma_1\big)\Phi_j^n.$ Plugging \eqref{prlefd1} into the semi-implicit 4cFD \eqref{eq:4cFD}, we obtain that the corresponding amplification factor $\eta_l$ satisfies
\begin{equation}
\label{prlefdX3}
\eta_{l}^{2}-2 i \tau \theta_{l} \eta_{l}-1=0, \quad l\in\mathcal{T}_N,
\end{equation}
with $\theta_{l}=-V_{\max}\pm A_{1,\max}$. The stability condition is equivalent to
\begin{equation}\label{prlefdX2}
\left|\eta_{l}\right| \leq 1 \Longleftrightarrow\left|\tau\theta_{l}\right| \leq 1, \quad l\in\mathcal{T}_N,
\end{equation}
which leads to
\[0<\tau\le\frac{1}{V_{\max}+A_{1,\max}}.\]
\hfill $\square$ \bigskip

\section{Error estimates for the 4cFD methods}
\label{sec:3}
In order to establish the error bounds for these two 4cFD methods, we assume the exact solution $\Phi$ of the Dirac equation \eqref{eq:Dirac_1D} satisfies
\begin{equation*}
\begin{split}
\Phi \in C^3([0, T_0];(L^{\infty}(\Omega))^2) & \cap C^2([0, T_0];(W^{1, \infty}_p(\Omega))^2)  \cap C^1([0, T_0];(W^{2, \infty}_p(\Omega))^2) \\
& \cap C([0, T_0];(W^{5, \infty}_p(\Omega))^2),
\end{split}
\end{equation*}
and
\begin{equation*}
(A)  \quad \left\|\frac{\partial^{r+s}}{\partial t^r\partial x^s}\Phi\right\|_{L^{\infty}([0, T_0];(L^{\infty}(\Omega))^2)} \lesssim \frac{1}{\varepsilon^{r}}, \quad 0 \leq r \leq 3,\, 0 \leq r+s \leq 5,\, 0 < \varepsilon \leq 1,
\end{equation*}
where $W^{m, \infty}_p = \{u | u \in W^{m, \infty}(\Omega),\ \partial^l_x u(a)= \partial^l_x u(b), \ l=0, \cdots, m-1\}$ for $m \geq 1$ and the boundary values are understood in the trace sense. In addition, we assume the electric potential $V \in C(\overline{\Omega}_{T})$ and magnetic potential $A_1 \in C(\overline{\Omega}_{T})$, and denote
\begin{equation*}
(B) \quad V_{\max}:= \max_{(t, x) \in \overline{\Omega}_{T}}|V(t, x)|,\quad A_{1, \max}:= \max_{(t, x) \in \overline{\Omega}_{T}}|A_1(t, x)|.
\end{equation*}

\subsection{The main results}
Let $\Phi^n_j$ be the numerical approximations of the 4cFD methods and define the grid error function $\mathbf{e}^n = (\mathbf{e}^n_0, \mathbf{e}^n_1,\cdots,\mathbf{e}^n_N)^T \in X_N$ as:
\begin{equation}\label{eq:e}
\mathbf{e}^n_j : = \Phi(t_n, x_j) - \Phi^n_j, \quad j=0,1, \cdots, N, \quad n\geq 0,
\end{equation}
then we could establish the error estimates for the 4cFD methods in the following theorems.

\begin{theorem}
\label{thm:phi2}
Under the assumptions  (A) and (B), there exist constants $h_0 > 0$ and $\tau_0 > 0$ sufficiently small and independent of $\varepsilon$, such that for any $0 < \varepsilon \leq 1$, when $0 < h \leq h_0$ and $0 < \tau \leq \tau_0$, we have the error bound for the implicit 4cFD \eqref{eq:im4cFD} with \eqref{eq:bic}
\begin{equation}
\label{eq:im4cFD_error}
\|\mathbf{e}^n\|_{l^2} \lesssim \frac{h^4}{\varepsilon} + \frac{\tau^2}{\varepsilon^3}, \quad 0 \leq n \leq \frac{T_0}{\tau}.
\end{equation}
\end{theorem}

\begin{theorem}
\label{thm:phi}
Under the assumptions (A) and (B), there exist constants $h_0 > 0$ and $\tau_0 > 0$ sufficiently small and independent of $\varepsilon$, such that for any $0 < \varepsilon \leq 1$, when $0 < h \leq h_0$ and $0 < \tau \leq \tau_0$ and under the stability condition \eqref{eq:stability}, we have the error bound for the semi-implicit 4cFD \eqref{eq:4cFD} with \eqref{eq:bic} and \eqref{eq:phi1}
\begin{equation}
\label{eq:4cFD_error}
\|\mathbf{e}^n\|_{l^2} \lesssim \frac{h^4}{\varepsilon} + \frac{\tau^2}{\varepsilon^3}, \quad 0 \leq n \leq \frac{T_0}{\tau}.
\end{equation}
\end{theorem}

Additionally, for the physical observables including the total density and current density, we also have the following error bounds.

 \begin{corollary}
\label{cor:rhoJ}
Under the assumptions (A) and (B), with the initial and boundary conditions \eqref{eq:bic}, \eqref{eq:phi1} and the corresponding stability condition for the implicit and semi-implicit 4cFD methods, there exist constants $h_0 > 0$ and $\tau_0 > 0$ sufficiently small and independent of $\varepsilon$, such that for any $0 < \varepsilon \leq 1$, when $0 < h \leq h_0$ and $0 < \tau \leq \tau_0$, we have the following error estimates on the total density and current density
\begin{equation*}
\|\rho^n - \rho(t_n, \cdot)\|_{l^2} \lesssim \frac{h^4}{\varepsilon} + \frac{\tau^2}{\varepsilon^3}, \quad \|{\bf J}^n - {\bf J}(t_n, \cdot)\|_{l^2} \lesssim \frac{h^4}{\varepsilon} + \frac{\tau^2}{\varepsilon^3}, \quad 0 \leq n \leq \frac{T_0}{\tau},
\end{equation*}
where $\rho^n$ and ${\bf J}^n$ are obtained from the numerical solution $\Phi^n$ via \eqref{eq:rho}, \eqref{eq:J} with $d=1$.
\end{corollary}

According to Theorems \ref{thm:phi2}, \ref{thm:phi} and Corollary \ref{cor:rhoJ}, we could obtain the spatial/temporal resolution capacity of the 4cFD methods for the Dirac equation \eqref{eq:Dirac_21} in the massless and nonrelativistic regime. Given an accuracy bound  $\delta_0 > 0$, the $\varepsilon$-scalability of the 4cFD methods is

\begin{equation*}
h = O(\varepsilon^{1/4}\sqrt{\delta_0}) = O(\varepsilon^{1/4}),\quad \tau = O(\varepsilon^{3/2}\sqrt{\delta_0}) = O(\varepsilon^{3/2}), \quad 0 < \varepsilon \leq 1.
\end{equation*}
Compared with the commonly used standard finite difference methods \cite{BCJT,MY,MY2}, the 4cFD methods could improve the spatial resolution capacity or attain higher order accuracy in space for a given mesh size, in other words, they need less grid points while maintaining the same accuracy.

\subsection{Proof of Theorem \ref{thm:phi2}}
\emph{Proof.}
Define the local truncation error ${\xi}^n=({\xi}_0^n,{\xi}_1^n,\ldots,
{\xi}_N^n)^T\in X_N$ of the implicit 4cFD \eqref{eq:im4cFD} with \eqref{eq:bic} for $n\ge0$ and $0\le j\le N-1,$
\begin{align}
\label{ccff1}
{\xi}_{j}^{n} := & \ i \delta_{t}^+ \Phi(t_n,x_j)+ \frac{i}{\varepsilon} \sigma_{1} \mathcal{A}^{-1}_h\delta_{x}\frac{\Phi(t_{n},x_j)+\Phi(t_{n+1},x_j)}{2} \nonumber \\
& \ -\Big(\frac{\sigma_{3}}{\varepsilon} +V_{j}^{n+\frac{1}{2}} I_{2}-A_{1, j}^{n+\frac{1}{2}} \sigma_{1}\Big) \frac{\Phi(t_{n},x_j)+\Phi(t_{n+1},x_j)}{2}.
\end{align}
Under the assumptions (A) and (B), by using the Taylor expansion, we could obtain for $n \geq 0$,
\begin{align}
|{\xi}_{j}^{n}| \le \ & \frac{\tau^2}{6} \|\partial_{ttt} \Phi\|_{l^\infty} + \frac{h^4}{120\varepsilon}\| \partial_{xxxxx}\Phi\|_{l^\infty}+ \frac{\tau^2}{4\varepsilon} \|\partial_{xtt}\Phi\|_{l^\infty}  \nonumber \\
&  + \frac{\tau^2}{4}\left(\frac{1}{\varepsilon} + V_{max} + A_{1,max} \right)\| \partial_{tt}\Phi\|_{l^\infty},
\end{align}
which leads to the following bounds
\begin{equation}\label{ccff4}
\left\|{\xi}^{n}\right\|_{l^{\infty}} \lesssim \frac{h^4}{\varepsilon}+\frac{\tau^2}{\varepsilon^3},\quad
\left\|{\xi}^{n}\right\|_{l^{2}} \lesssim\left\|{\xi}^{n}\right\|_{l^{\infty}} \lesssim \frac{h^{4}}{\varepsilon}+\frac{\tau^2}{\varepsilon^3}.
\end{equation}
Subtracting \eqref{eq:im4cFD} from \eqref{ccff1} and noticing \eqref{eq:e}, we obtain the error function with $n\ge0$ and $0\le j\le N-1$
\begin{equation}\label{ccff41}
 i \delta_{t}^+ \mathbf{e}_{j}^{n}=-\frac{i}{\varepsilon} \sigma_{1} \mathcal{A}^{-1}_h \delta_{x} \mathbf{e}_{j}^{n+\frac{1}{2}}+
\left(\frac{\sigma_{3}}{\varepsilon}+V_{j}^{n+\frac{1}{2}} I_{2}-A_{1, j}^{n+\frac{1}{2}} \sigma_{1}\right) \mathbf{e}_{j}^{n+\frac{1}{2}}+{\xi}_{j}^{n},
\end{equation}
where the boundary and initial conditions are given as
\begin{equation}\label{bce}
\mathbf{e}^n_0 = \mathbf{e}^n_N,\quad \mathbf{e}^n_{-1} = \mathbf{e}^n_{N-1},\quad n \geq 0, \quad \mathbf{e}^0_j = {\bf 0}, \quad j = 0, 1, \cdots, N.
\end{equation}
Noticing \eqref{ccff4}, multiplying $h\tau(\mathbf{e}_j^{n+1}+\mathbf{e}_j^{n})^*$ from the left on both sides of \eqref{ccff41}, summing up from $j = 0$ to $j = N-1$, and taking the imaginary part, we arrive at
\begin{align}
\left\|\mathbf{e}^{n+1}\right\|_{l^{2}}^{2}-\left\|\mathbf{e}^{n}\right\|_{l^{2}}^{2}
 & =\tau \operatorname{Im}\bigg[h \sum_{j=0}^{N-1}\left(\mathbf{e}_{j}^{n+1}+\mathbf{e}_{j}^{n}\right)^{*} {{\xi}}_{j}^{n}\bigg] \nonumber\\
 & \lesssim \tau\left(\|\mathbf{e}^{n+1}\|_{l^2}^2+\|\mathbf{e}^{n}\|_{l^2}^2\right)
 +\tau\left(\frac{h^4}{\varepsilon} + \frac{\tau^2}{\varepsilon^3}\right)^2, \quad n \geq 0.
 \end{align}
Summing up the above equalities for $n=0,1,2,\cdots,m-1,$ we obtain
\begin{equation}\label{ccff9}
\left\|\mathbf{e}^{m}\right\|_{l^{2}}^{2}-\left\|\mathbf{e}^{0}\right\|_{l^{2}}^{2}
\lesssim\tau \sum_{s=0}^{m}\|\mathbf{e}^{s}\|_{l^2}^2
+m\tau\left(\frac{h^{4}}{\varepsilon}+\frac{\tau^2}{\varepsilon^3}\right)^2,
\quad 1 \leq m \leq \frac{T_0}{\tau},
\end{equation}
where $\left\|\mathbf{e}^{0}\right\|_{l^{2}}^{2}=0.$ By taking $\tau_0$ sufficiently small and using the discrete Gronwall's inequality, we get
\begin{equation}\label{ccff10}
\left\|\mathbf{e}^{m}\right\|^2_{l^{2}}
\lesssim \left(\frac{h^{4}}{\varepsilon}+\frac{\tau^{2}}{\varepsilon^{3}}\right)^2, \quad 1 \leq m \leq \frac{T_0}{\tau},
\end{equation}
which completes the proof of the error bound \eqref{eq:im4cFD_error}.
\hfill $\square$ \bigskip

\subsection{Proof of Theorem \ref{thm:phi}}
\emph{Proof.}
Define the local truncation error ${\xi}^n = ({\xi}^n_0,{\xi}^n_1, \cdots, {\xi}^n_N)^T \in X_N$ of the semi-implicit 4cFD \eqref{eq:4cFD} with \eqref{eq:bic} and \eqref{eq:phi1} for $n \geq 1$ and $0 \leq j \leq N-1,$
\begin{align}
\label{eq:xi}	
{\xi}^n_j := &\ i\delta_t\Phi(t_n, x_j) + \left( \frac{i}{\varepsilon} \sigma_1\mathcal{A}^{-1}_h\delta_x - \frac{1}{\varepsilon}\sigma_3\right) \frac{\Phi(t_{n-1}, x_j) + \Phi(t_{n+1}, x_j)}{2} \nonumber\\
 &\ - \left( V^{n}_j I_2 - A^{n}_{1, j} \sigma_1  \right)\Phi(t_n, x_j),
\end{align}
and
\begin{align}
\label{eq:xi0}	
{\xi}^0_j := \ & i\delta_t^+ \Phi_0(x_j) + \frac{i}{\varepsilon}\sigma_ 1\Phi'_0(x_j)-\left(\frac{1}{\varepsilon}\sigma_3 +V^{0}_j I_2 - A^{0}_{1, j} \sigma_1 \right)\Phi_0(x_j).
\end{align}
Noticing \eqref{eq:Dirac_1D} and the assumptions $(A)$ and $(B)$, and applying the Taylor expansions \eqref{eq:xi}, \eqref{eq:xi0}, we have for $n \geq 1$ and $0\le j \le N-1$,
\begin{equation*}
\begin{split}
|{\xi}^n_j |\le&\  \frac{\tau^2}{6}\|\partial_{ttt}\Phi\|_{l^\infty} + \frac{\tau^2}{2\varepsilon} \|\partial_{xtt}\Phi\|_{l^\infty} + \frac{h^4}{120\varepsilon}\|\partial_{xxxxx}\Phi\|_{l^\infty}  +\frac{\tau^2}{2\varepsilon}\|\partial_{tt}\Phi\|_{l^\infty}, \\
|\xi^0_j| \le &\  \frac{\tau}{2} \|\partial_{tt}\Phi\|_{l^\infty},
\end{split}
\end{equation*}
which imply the following error bounds
\begin{align}
\label{eq:xib}
&\|{\xi}^0\|_{l^{\infty}}	= \max_{0 \leq j \leq N-1} |\xi^0_j| \lesssim \frac{\tau}{\varepsilon^2}, \quad \|{\xi}^0\|_{l^2} \lesssim \|{\xi}^0\|_{l^{\infty}} \lesssim \frac{\tau}{\varepsilon^2}, \nonumber \\
&\|{\xi}^n\|_{l^{\infty}}	= \max_{0 \leq j \leq N-1} |\xi^n_j| \lesssim \frac{h^4}{\varepsilon} + \frac{\tau^2}{\varepsilon^3}, \quad \|{\xi}^n\|_{l^2} \lesssim \|{\xi}^n\|_{l^{\infty}} \lesssim  \frac{h^4}{\varepsilon} + \frac{\tau^2}{\varepsilon^3}, \quad n \geq 1.
\end{align}
By the definition of the error function \eqref{eq:e}, we get for $0 \leq j \leq N-1$ and $n \geq 1$,
\begin{equation}
i\delta_t\mathbf{e}^n_j = \left(-\frac{i}{2\varepsilon}\sigma_1 \mathcal{A}^{-1}_h\delta_x +\frac{1}{2\varepsilon} \sigma_3\right)\left({\bf e}^{n+1}_j + {\bf e}^{n-1}_j\right) + \left( V^{n}_j I_2 - A^{n}_{1, j} \sigma_1\right)\mathbf{e}^{n}_j  + {\xi}^n_j,
\label{eq:error}
\end{equation}
and the initial and boundary conditions are same as \eqref{bce}. For the first step, we have $\frac{i}{\tau} {\bf e}^1_j = \xi^0_j$ for $j = 0, 1, \cdots, N$, then
\begin{equation}
\|{\bf e}^1\|_{l^2} = \tau \|\xi^0_j\|_{l^2} \lesssim \frac{\tau^2}{\varepsilon^2}	\lesssim \frac{\tau^2}{\varepsilon^3}.
\label{eq:e1}
\end{equation}
Denote
\begin{align}
\mathcal{E}^{n+1} = & \ \|{\bf e}^{n+1}\|^2_{l^2} + \frac{\tau h}{\varepsilon} \text{Re} \left(\sum^{N-1}_{j=0} ({\bf e}^{n+1}_j)^{\ast} \sigma_1 \mathcal{A}^{-1}_h \delta_x {\bf e}^{n+1}_j +\sum^{N-1}_{j=0}  ({\bf e}^{n}_j)^{\ast} \sigma_1 \mathcal{A}^{-1}_h \delta_x {\bf e}^{n}_j\right) \nonumber \\
&\  + \|{\bf e}^{n}\|^2_{l^2} -2\tau \text{Im} \left(h\sum^{N-1}_{j=0} ({\bf e}^{n+1}_j)^{\ast} \left(V^n_j I_2 -A^n_{1, j}\sigma_1\right){\bf e}^{n}_j\right), \quad n \geq 0.
\label{eq:matE}
\end{align}
Multiplying \eqref{eq:error} from the left on both sides with $2\tau h \left(\mathbf{e}^{n+1}_j + \mathbf{e}^{n-1}_j\right)^{\ast}$, summing up for $j = 0, 1, \cdots, N-1$, and taking the imaginary part,  then using Cauchy inequality and noticing \eqref{eq:xib}, we get
\begin{align}
\mathcal{E}^{n+1} - \mathcal{E}^n = 2\tau h \text{Im}\left[\sum^{N-1}_{j = 0}\left(\mathbf{e}^{n+1}_j + \mathbf{e}^{n-1}_j\right)^{\ast}{\xi}^n_j\right] \nonumber \\
\lesssim \tau \left(\mathcal{E}^{n+1} + \mathcal{E}^n\right) + \tau \left(\frac{h^4}{\varepsilon} + \frac{\tau^2}{\varepsilon^3}\right)^2.
\end{align}
Under the stability condition \eqref{eq:stability}, taking $\frac{\tau}{h\varepsilon}\leq \frac{1}{4}$ and $\tau \left(V_{\max} + A_{1, \max}\right) \leq \frac{1}{4}$, we could derive
\begin{equation}
\frac{1}{2}\left(\|{\bf e}^{n+1}\|^2_{l^2} +\|{\bf e}^{n}\|^2_{l^2}\right) \leq \mathcal{E}^{n+1}	 \leq \frac{3}{2}\left(\|{\bf e}^{n+1}\|^2_{l^2} +\|{\bf e}^{n}\|^2_{l^2}\right), \quad n \geq 0,
\label{eq:de}
\end{equation}
by using Cauchy inequality. Combining \eqref{eq:e1} and the definition of $\mathcal{E}^n$ \eqref{eq:matE}, we arrive at
\begin{equation}
\mathcal{E}^1 \lesssim \left(\frac{h^4}{\varepsilon} + \frac{\tau^2}{\varepsilon^3}\right)^2.	
\end{equation}
Thus, if we take $\tau_0$ sufficiently small, the discrete Gronwall's inequality would yield
\begin{equation}
\mathcal{E}^m \lesssim  \left(\frac{h^4}{\varepsilon} + \frac{\tau^2}{\varepsilon^3}\right)^2, \quad 1 \leq m \leq \frac{T_0}{\tau},
\end{equation}
which implies the error bound \eqref{eq:4cFD_error} combining with the inequality \eqref{eq:de}.
\hfill $\square$ \bigskip

\section{Numerical results}
\label{sec:4}
In this section, we first show the example in 1D to confirm the accuracy and $\varepsilon$-scalability of the 4cFD methods for the Dirac equation \eqref{eq:Dirac_21} in the massless and nonrelativistic regime. Then we present the dynamics of the Dirac equation \eqref{eq:Dirac_21} in 2D with different potentials by the semi-implicit 4cFD method.

\subsection{Spatial and temporal resolution}
In this subsection, we take $\Omega = (0, 2\pi)$ and the electromagnetic potentials as
\begin{equation}\label{pot1}
V(t, x) = \frac{1}{1+\sin^2(x)},\quad A_1(t,x) = \sin(2x), \quad x \in \Omega, \quad t \geq 0,
\end{equation}
and choose the following initial data
\begin{equation}
\phi_1(0, x) = \frac{1}{1+\sin^2(x)},\quad \phi_2(0, x) = \frac{1}{3+\cos(x)},\quad x \in \Omega.
\end{equation}
Since the exact solution is unknown, we use the time-splitting Fourier pseudospectral \cite{BCJT} method  with a fine mesh size $h_e = \pi/512$ and a very small time step $\tau_e = 10^{-6}$ to obtain the `reference' solution numerically. Denote $\Phi^n_{h, \tau}$ as the numerical solution obtained by the 4cFD methods with the mesh size $h$ and time step $\tau$. In order to quantify the numerical errors of the 4cFD schemes, we use the relative errors for the wave function $\Phi$, the total density $\rho$ and current density $\bf J$
\begin{equation}
e^{\Phi}_{h, \tau}(t_n) = \frac{\|\Phi^n_{h, \tau} - \Phi(t_n, \cdot)\|_{l^2}}{\|\Phi(t_n, \cdot)\|_{l^2}},\quad
 e^{\rho}_{h, \tau}(t_n) = \frac{\|\rho^n_{h, \tau} - \rho(t_n, \cdot)\|_{l^2}}{\|\rho(t_n, \cdot)\|_{l^2}},\quad e^{\bf J}_{h, \tau}(t_n) = \frac{\|{\bf J}^n_{h, \tau} - {\bf J}(t_n, \cdot)\|_{l^2}}{{\|\bf J}(t_n, \cdot)\|_{l^2}}.
\end{equation}


\begin{table}[ht!]
\renewcommand{\arraystretch}{1.5}
\caption{Spatial errors for the wave function $e^{\Phi}_{h, \tau}(t =2)$ by using the 4cFD method}
\vspace{5pt}
\centering
{\normalsize
\begin{tabular}{ccccc}
\hline
$e^{\Phi}_{h, \tau}(t=2)$ &$h_0 = \frac{\pi}{16}$ & $h_0/2 $ &$h_0/2^2 $ & $h_0/2^3$  \\
\midrule[0.5pt]
$\varepsilon_0 = 1$ & 4.68E-3 & \bf{2.40E-4} & 1.44E-5 & 8.94E-7  \\
order & - & \bf{4.29} & 4.06 & 4.01 \\
$\varepsilon_0 / 2^4$  & 1.57E-2 & 1.24E-3 & \bf{7.62E-5} & 4.72E-6 \\
order & - & 3.66 & \bf{4.02} & 4.01 \\
$\varepsilon_0 / 2^8$  & 6.18E-2 & 1.47E-2 & 1.19E-3 & \bf{7.49E-5} \\
order & - & 2.06 & 3.63 & \bf{3.99} \\
\hline
\end{tabular}}
\label{tab:phi_h}
\end{table}

\begin{table}[ht!]
\renewcommand{\arraystretch}{1.5}
\caption{Spatial errors for the total density $e^{\rho}_{h, \tau}(t =2)$ by using the 4cFD method}
\vspace{5pt}
\centering
{\normalsize
\begin{tabular}{ccccc}
\hline
$e^{\rho}_{h, \tau}(t=2)$ &$h_0 = \frac{\pi}{16}$ & $h_0/2 $ &$h_0/2^2 $ & $h_0/2^3$  \\
\midrule[0.5pt]
$\varepsilon_0 = 1$ & 4.31E-3 & \bf{2.06E-4} & 1.25E-5 & 7.73E-7  \\
order & - & \bf{4.39} & 4.04 & 4.02 \\
$\varepsilon_0 / 2^4$  & 1.58E-2 & 1.02E-3 & \bf{6.57E-5} & 4.08E-6 \\
order & - & 3.95 & \bf{3.96} & 4.01 \\
$\varepsilon_0 / 2^8$  & 8.59E-2 & 1.51E-2 & 1.49E-3 & \bf{9.30E-5} \\
order & - & 2.51 & 3.34 & \bf{4.00} \\
\hline
\end{tabular}}
\label{tab:rho_h}
\end{table}

\begin{table}[ht!]
\renewcommand{\arraystretch}{1.5}
\caption{Spatial errors for the current density $e^{\bf J}_{h, \tau}(t =2)$ by using the 4cFD method}
\vspace{5pt}
\centering
{\normalsize
\begin{tabular}{ccccc}
\hline
$e^{\bf J}_{h, \tau}(t=2)$ &$h_0 = \frac{\pi}{16}$ & $h_0/2 $ &$h_0/2^2 $ & $h_0/2^3$  \\
\midrule[0.5pt]
$\varepsilon_0 = 1$ & 8.09E-3 & \bf{3.74E-4} & 2.24E-5 & 1.39E-6  \\
order & - & \bf{4.44} & 4.06 & 4.01 \\
$\varepsilon_0 / 2^4$  & 3.58E-2 & 3.14E-3 & \bf{1.99E-4} & 1.24E-5 \\
order & - & 3.51 & \bf{3.98} & 4.00 \\
$\varepsilon_0 / 2^8$  & 9.75E-2 & 3.12E-2 & 2.13E-3 & \bf{1.36E-4} \\
order & - & 1.64 & 3.87 & \bf{3.97} \\
\hline
\end{tabular}}
\label{tab:J_h}
\end{table}

\begin{table}[htp!]
\renewcommand{\arraystretch}{1.5}
\caption{Temporal errors for the wave function $e^{\Phi}_{h, \tau}(t =2)$ by using the 4cFD method}
\vspace{5pt}
\centering
{\normalsize
\setlength{\tabcolsep}{1mm}
\begin{tabular}{cccccccc}
\hline
$e^{\Phi}_{h, \tau}(t=2)$ &$\tau_0 = 0.05$ & $\tau_0/2 $ & $\tau_0/2^2$ & $\tau_0/2^3$ & $\tau_0/2^4$ & $\tau_0/2^5$  & $\tau_0/2^6$ \\
\midrule[0.5pt]
$\varepsilon_0 = 1$ & \bf{1.71E-2} & 4.31E-3 & 1.08E-3 & 2.70E-4  & 6.74E-5 & 1.69E-5 & 4.22E-6 \\
order & \bf{-} & 1.99 & 2.00 & 2.00 & 2.00 & 2.00 & 2.00  \\
$\varepsilon_0 / 2^{2/3}$ & 4.02E-2 & \bf{1.04E-2} & 2.61E-3 & 6.52E-4 & 1.63E-4  & 4.08E-5 & 1.02E-5 \\
order & - & \bf{1.95} & 1.99 & 2.00 & 2.00 & 2.00 & 2.00  \\
$\varepsilon_0 / 2^{4/3}$ & 1.15E-1 & 3.27E-2 & \bf{8.42E-3} & 2.12E-3 & 5.30E-4 & 1.32E-4    & 3.31E-5  \\
order & - & 1.81 & \bf{1.96} & 1.99 & 2.00 & 2.01  & 2.00 \\
$\varepsilon_0 / 2^{2}$ & 3.03E-1 & 1.03E-1 & 2.92E-2 & \bf{7.50E-3} & 1.88E-3 & 4.71E-4 & 1.18E-4  \\
order & - & 1.56 & 1.82 & \bf{1.96} & 2.00 & 2.00 & 2.00  \\
$\varepsilon_0 / 2^{8/3}$ & 6.59E-1 & 2.83E-1 & 9.68E-2 & 2.75E-2 & \bf{7.03E-3} & 1.76E-3  & 4.41E-4  \\
order & - & 1.22 & 1.55 & 1.82 & \bf{1.97} & 2.00 & 2.00 \\
$\varepsilon_0 / 2^{10/3}$ & 1.49 &  6.06E-1 & 2.73E-1 & 9.37E-2 & 2.66E-2 & \bf{6.81E-3} & 1.71E-3 \\
order & - & 1.30 & 1.15 & 1.54 & 1.82 & \bf{1.97} & 1.99 \\
\hline
\end{tabular}}
\label{tab:phi_t}
\end{table}

\begin{figure}[ht!]
\begin{minipage}{0.5\textwidth}
\centerline{\includegraphics[width=6.5cm,height=5cm]{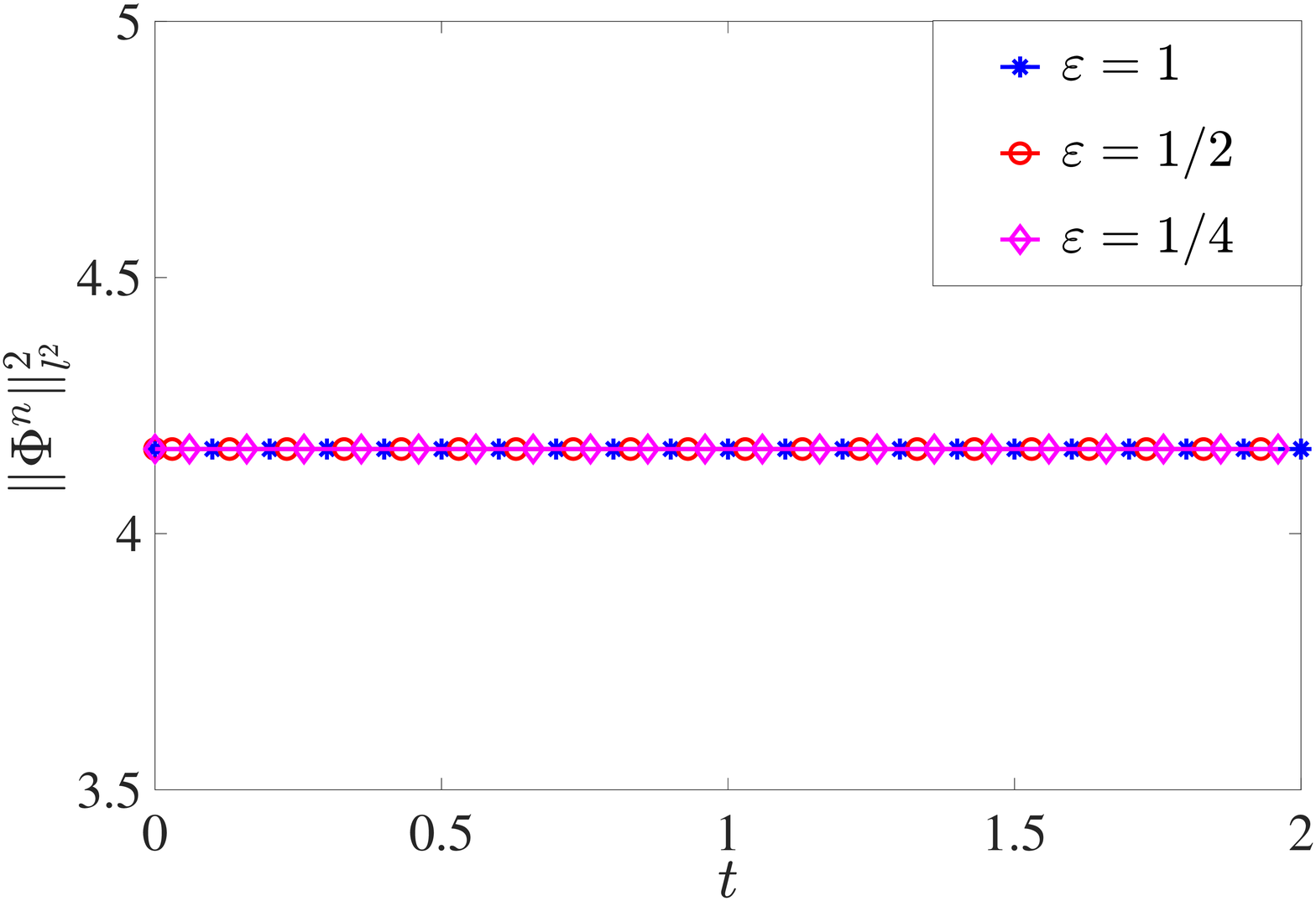}}
\end{minipage}
\begin{minipage}{0.48\textwidth}
\centerline{\includegraphics[width=6.5cm,height=5cm]{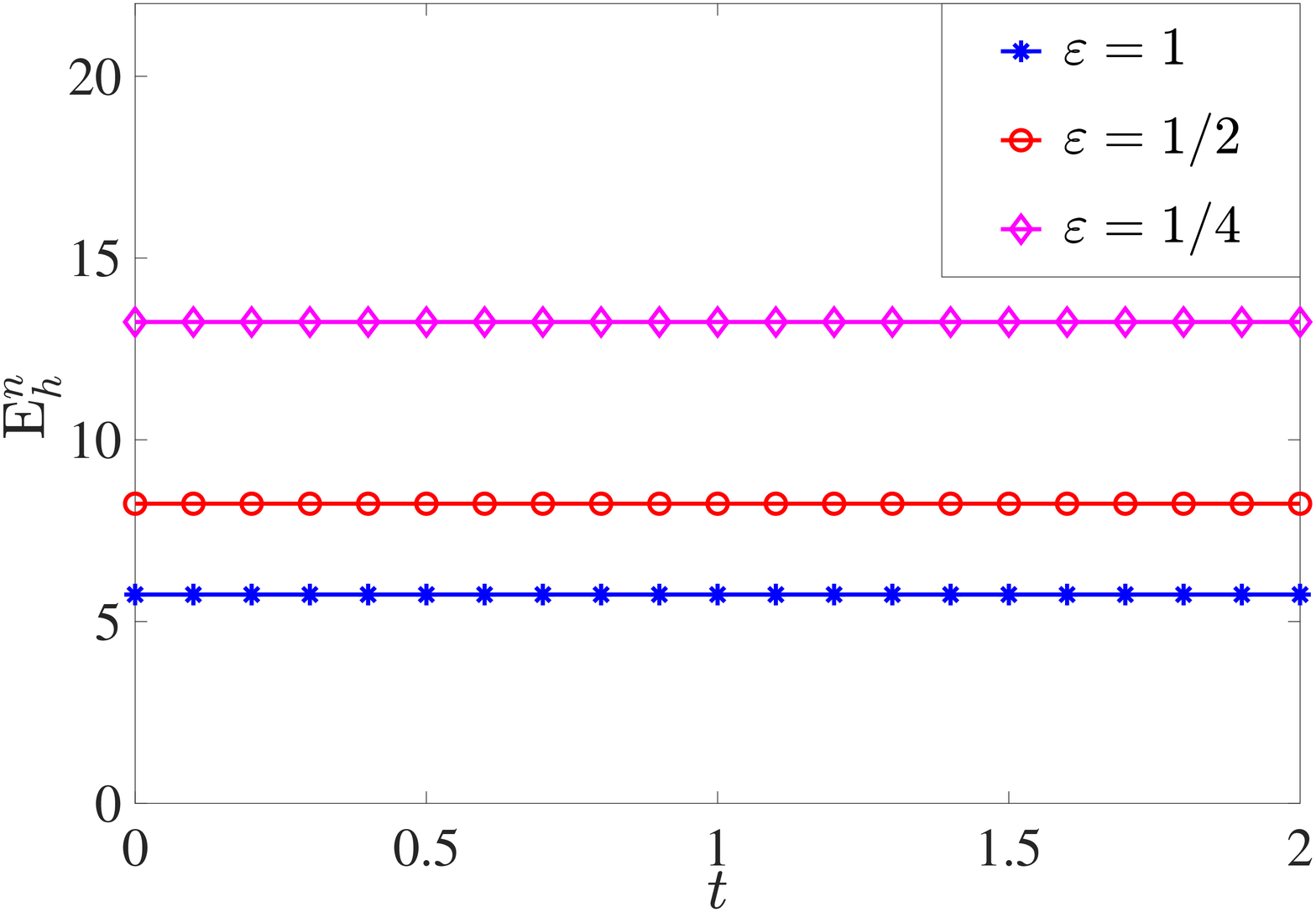}}
\end{minipage}
\caption{Evolution of the mass $\|\Phi^n\|^2_{l^2}$ (left) and energy $E^n_h$ (right) for the implicit 4cFD method with different $\varepsilon$}
\label{fig:ME}
\end{figure}
Tables \ref{tab:phi_h}--\ref{tab:J_h} show the relative spatial errors for the wave function $e^{\Phi}_{h, \tau}(t = 2)$, the total density $e^{\rho}_{h, \tau}(t = 2)$ and current density $e^{\bf J}_{h, \tau}(t = 2)$, respectively, by using the semi-implicit 4cFD method. Table \ref{tab:phi_t} shows the relative temporal errors for the wave function $e^{\Phi}_{h, \tau}(t = 2)$. The numerical results for the implicit 4cFD scheme are similar to the semi-implicit scheme and we omit the results here for brevity. From Tables \ref{tab:phi_h}--\ref{tab:phi_t}, we could observe that the 4cFD schemes are fourth-order accurate in space and second-order accurate in time (cf. the first row in Tables \ref{tab:phi_h}--\ref{tab:phi_t}). The $\varepsilon$-scalability of the 4cFD methods is $h = O(\varepsilon^{1/4})$ and $\tau = O(\varepsilon^{3/2})$, which is verified through the upper part of each table above the bold diagonal line. For the relative temporal errors of the total density and current density, we could observe similar results with the relative errors for the wave function. Fig. \ref{fig:ME} depicts the evolution of the mass and energy for the implicit 4cFD method and demonstrates that the implicit scheme preserves the mass and energy conservation and the mass is independent of the small parameter $\varepsilon$ while the energy becomes larger as $\varepsilon$ is smaller.
In summary, the numerical results confirm the error bounds for the wave function in Theorems \ref{thm:phi2}, \ref{thm:phi} and for the total density and current density in Corollary \ref{cor:rhoJ} for the Dirac equation in the massless and nonrelativistic regime.

\subsection{Dynamics of the Dirac equation in 2D}
In this subsection, we apply the semi-implicit 4cFD method to study the dynamics of the Dirac equation \eqref{eq:Dirac_21} in 2D with different electromagnetic potentials.

\medskip

\begin{figure}[h!]
\begin{minipage}{0.31\textwidth}
\centerline{\includegraphics[width=4.3cm,height=3.6cm]{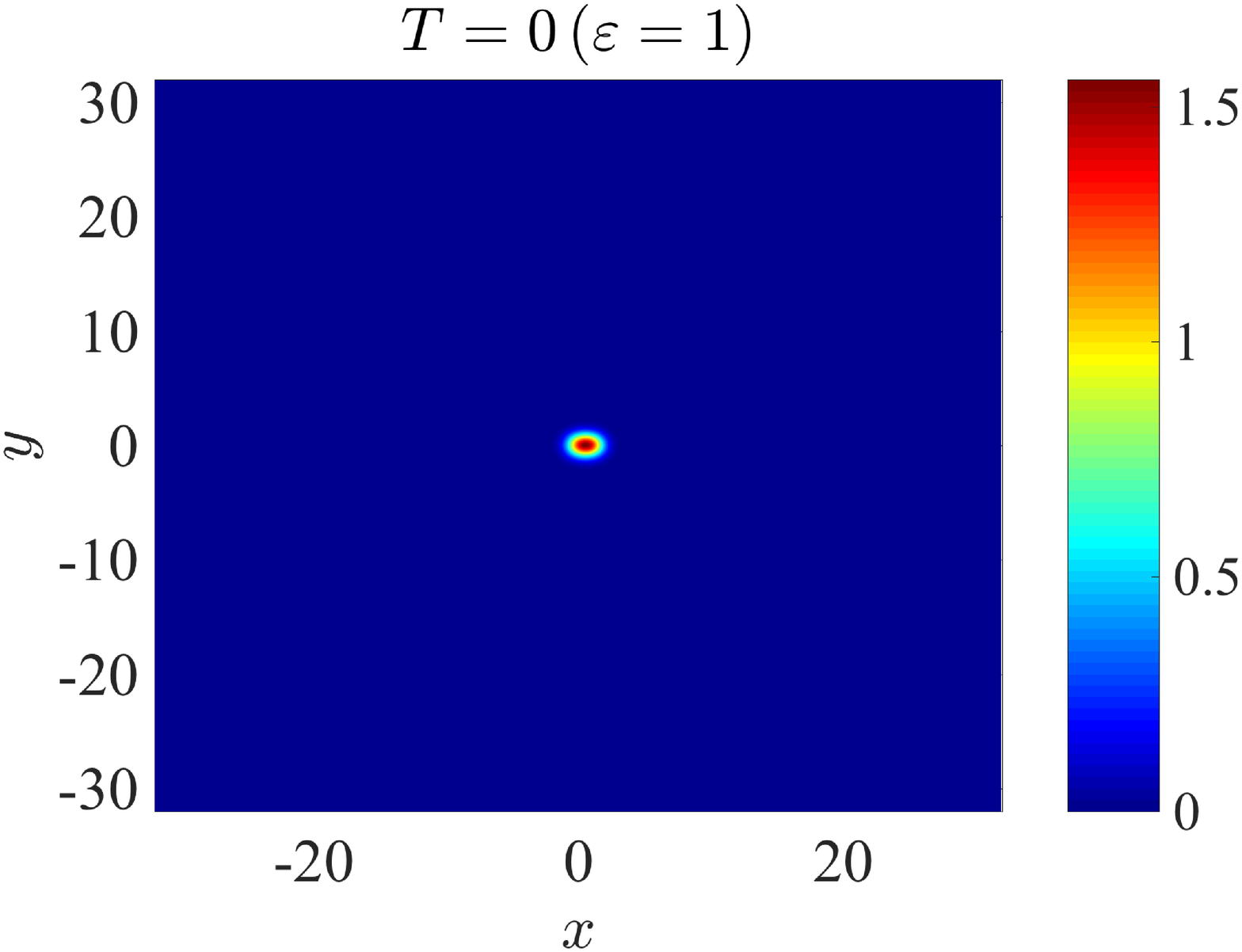}}
\end{minipage}
\begin{minipage}{0.31\textwidth}
\centerline{\includegraphics[width=4.3cm,height=3.6cm]{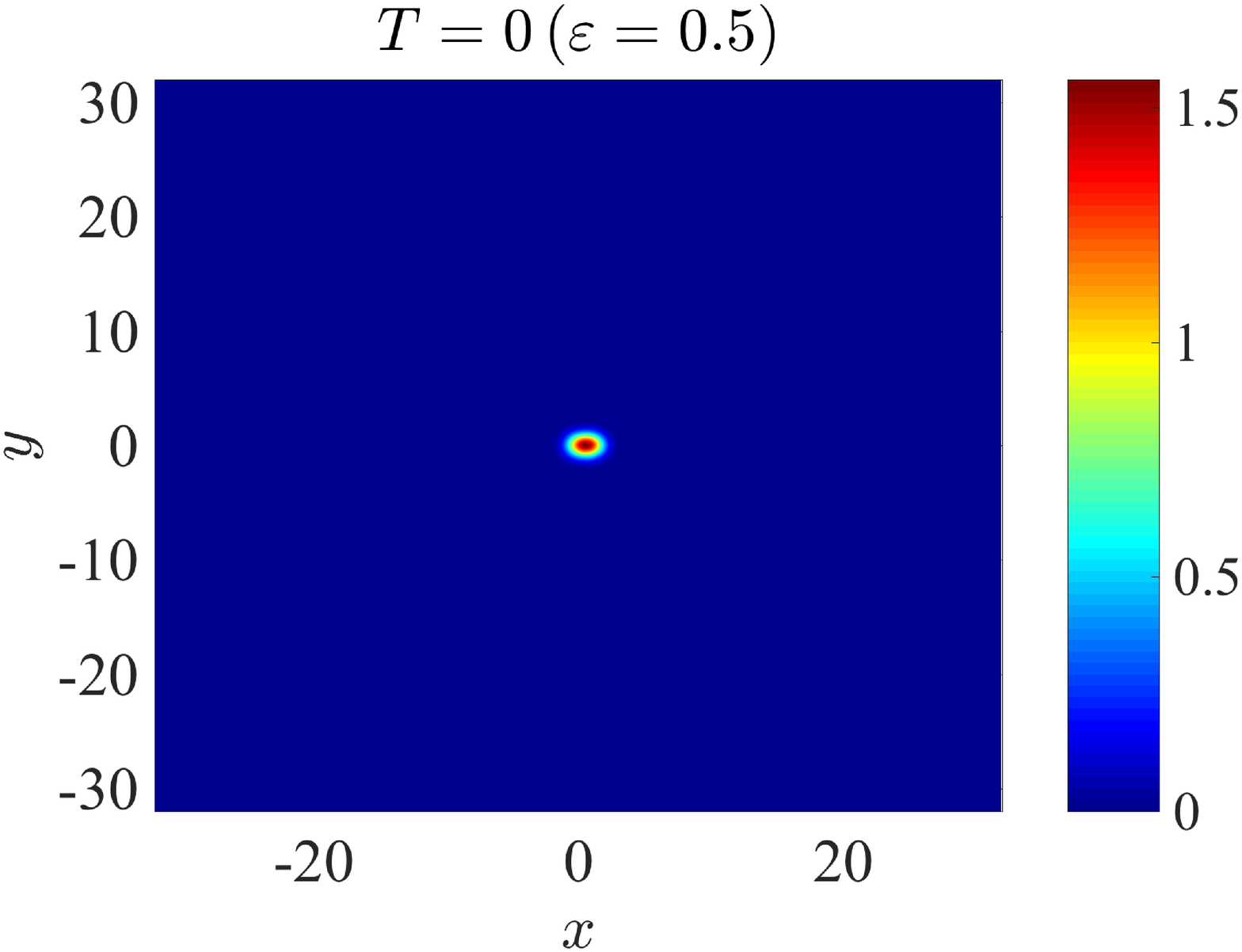}}
\end{minipage}
\begin{minipage}{0.31\textwidth}
\centerline{\includegraphics[width=4.3cm,height=3.6cm]{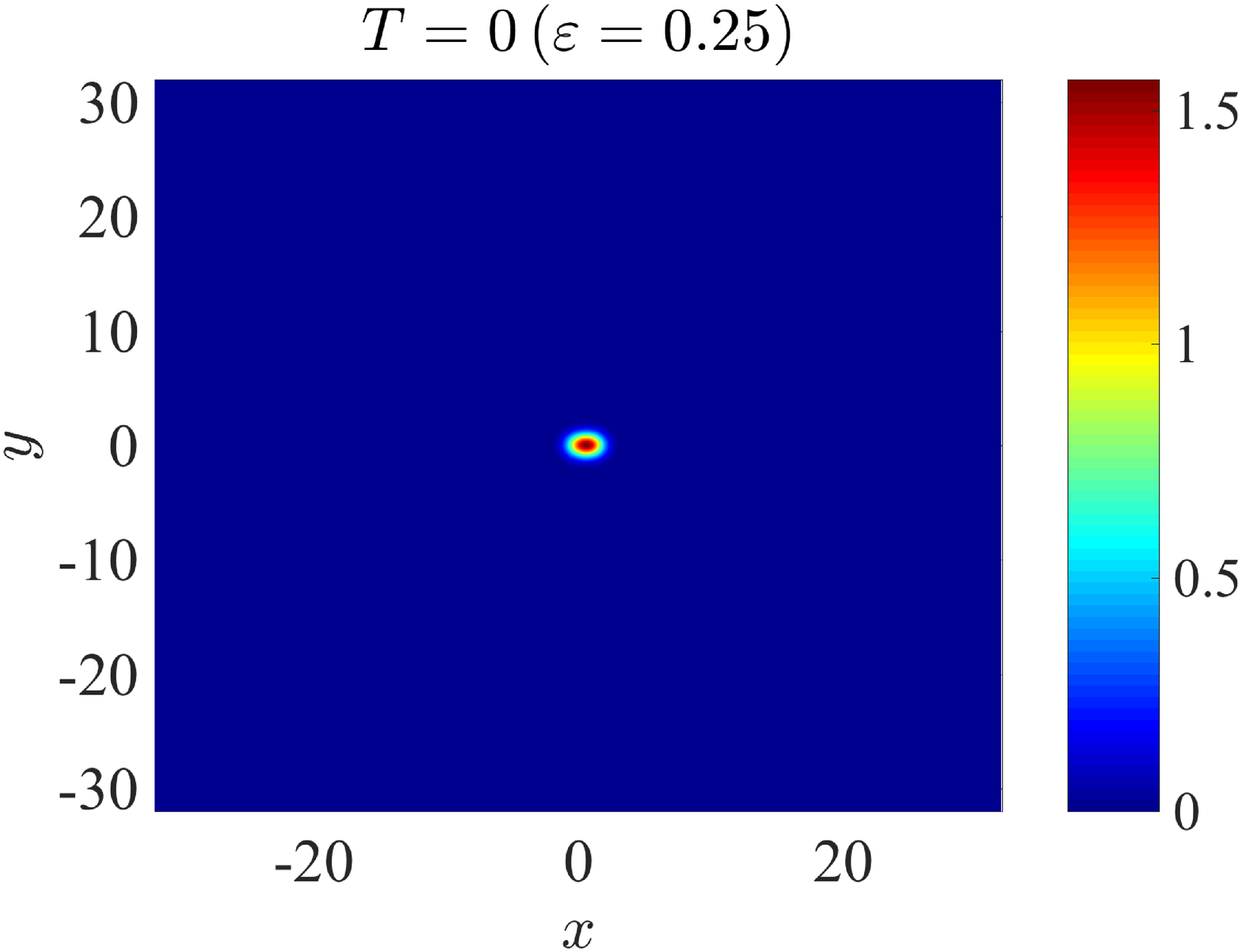}}
\end{minipage} \\
\begin{minipage}{0.31\textwidth}
\centerline{\includegraphics[width=4.3cm,height=3.6cm]{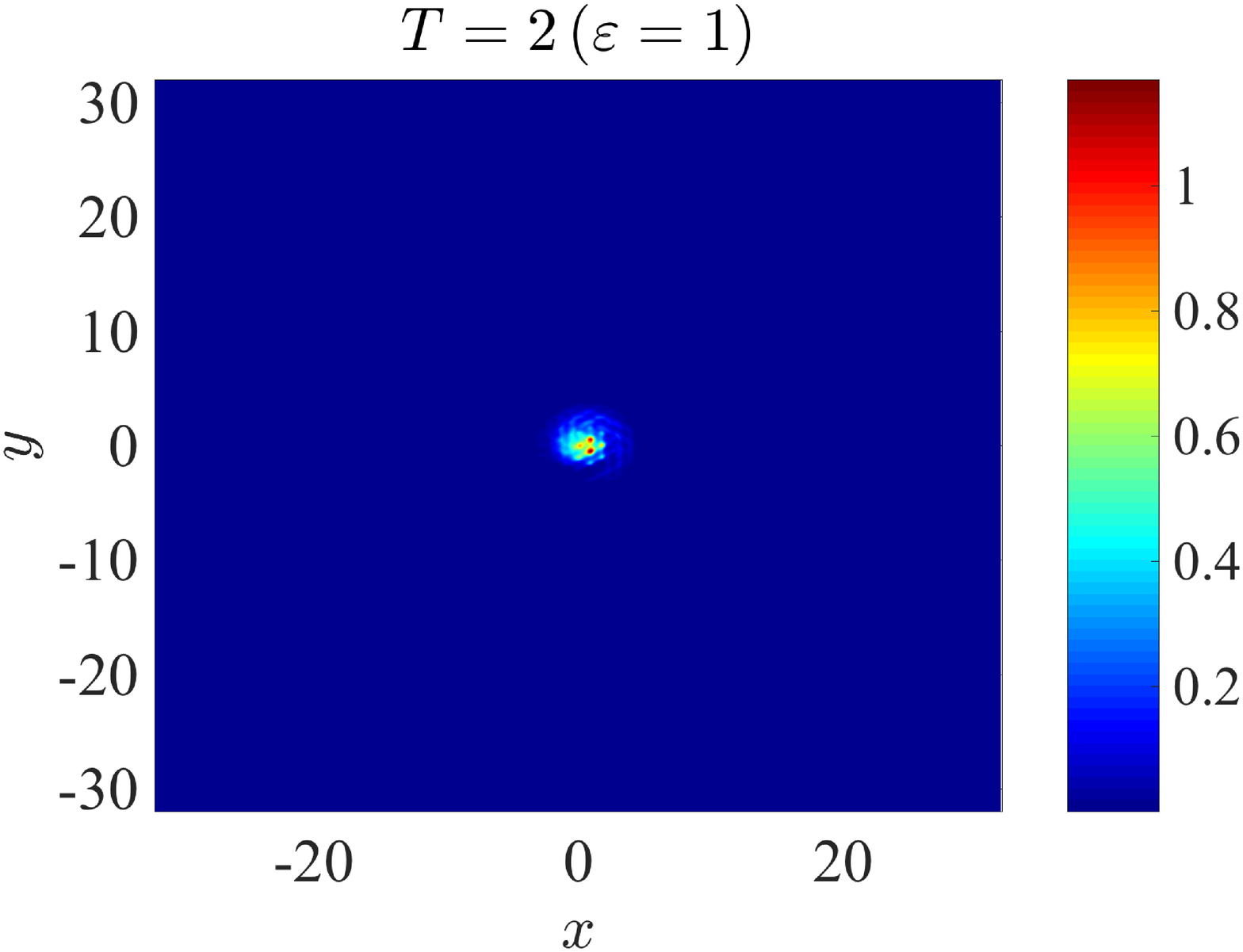}}
\end{minipage}
\begin{minipage}{0.31\textwidth}
\centerline{\includegraphics[width=4.3cm,height=3.6cm]{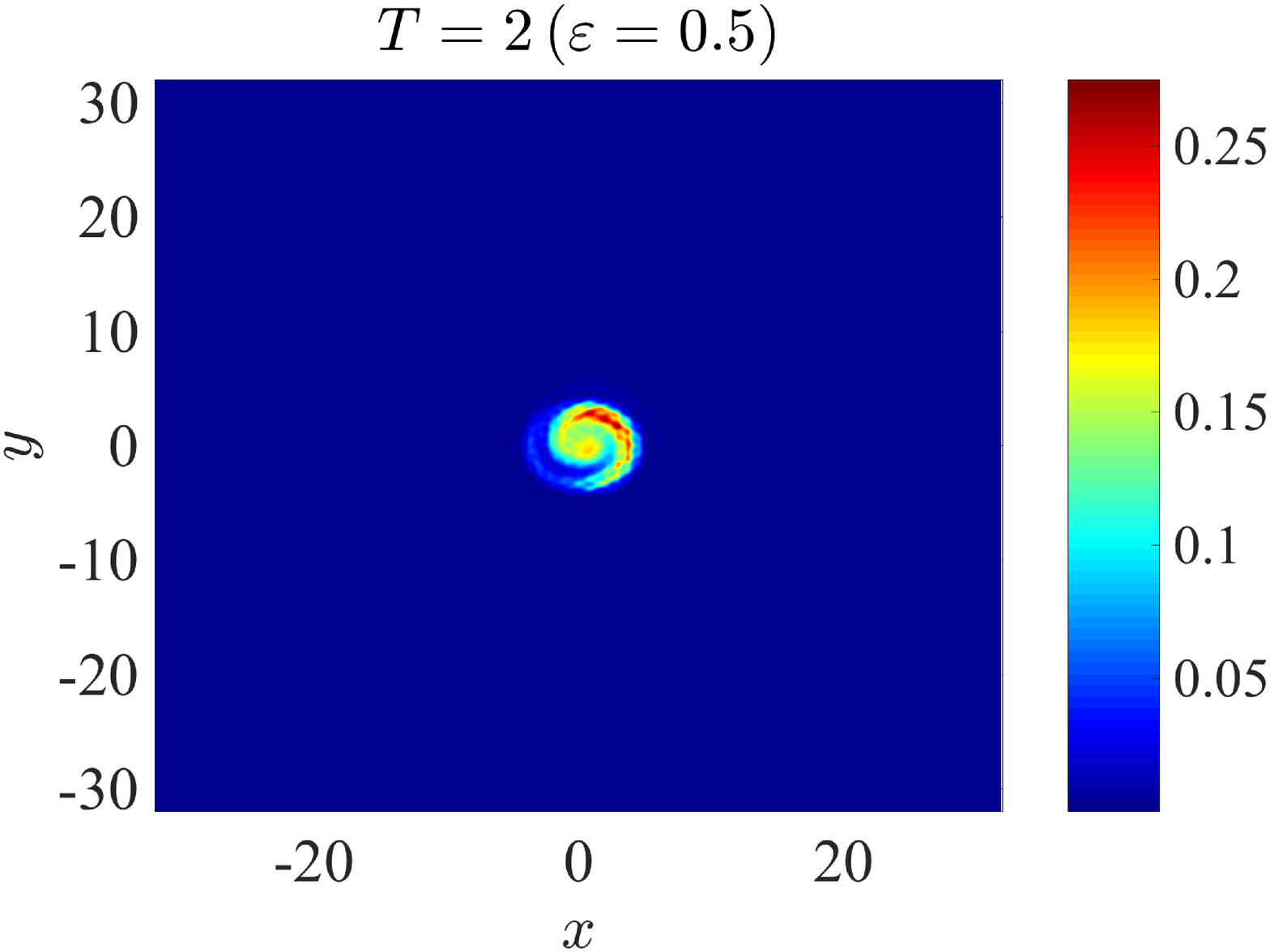}}
\end{minipage}
\begin{minipage}{0.31\textwidth}
\centerline{\includegraphics[width=4.3cm,height=3.6cm]{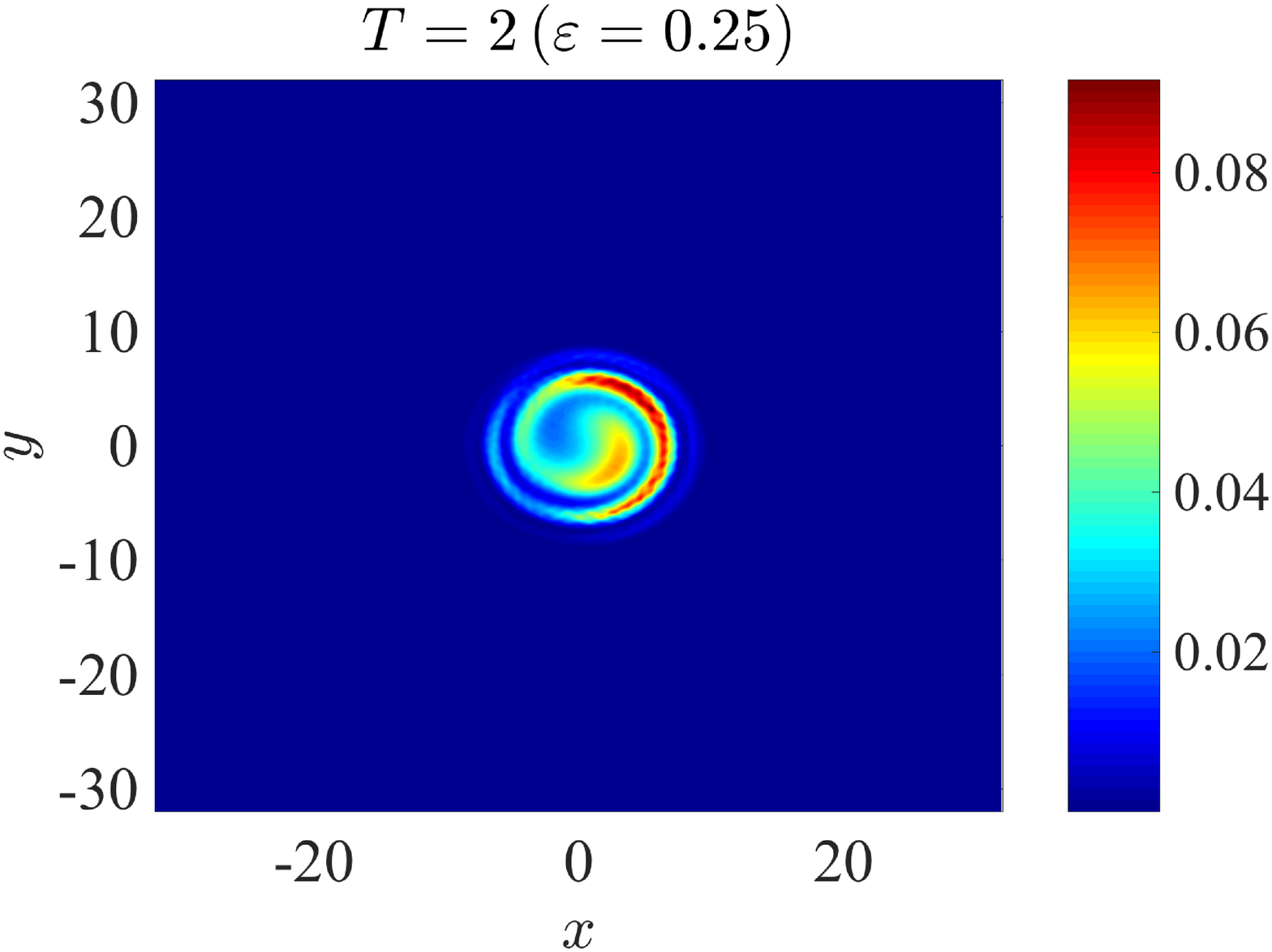}}
\end{minipage}\\
\begin{minipage}{0.31\textwidth}
\centerline{\includegraphics[width=4.3cm,height=3.6cm]{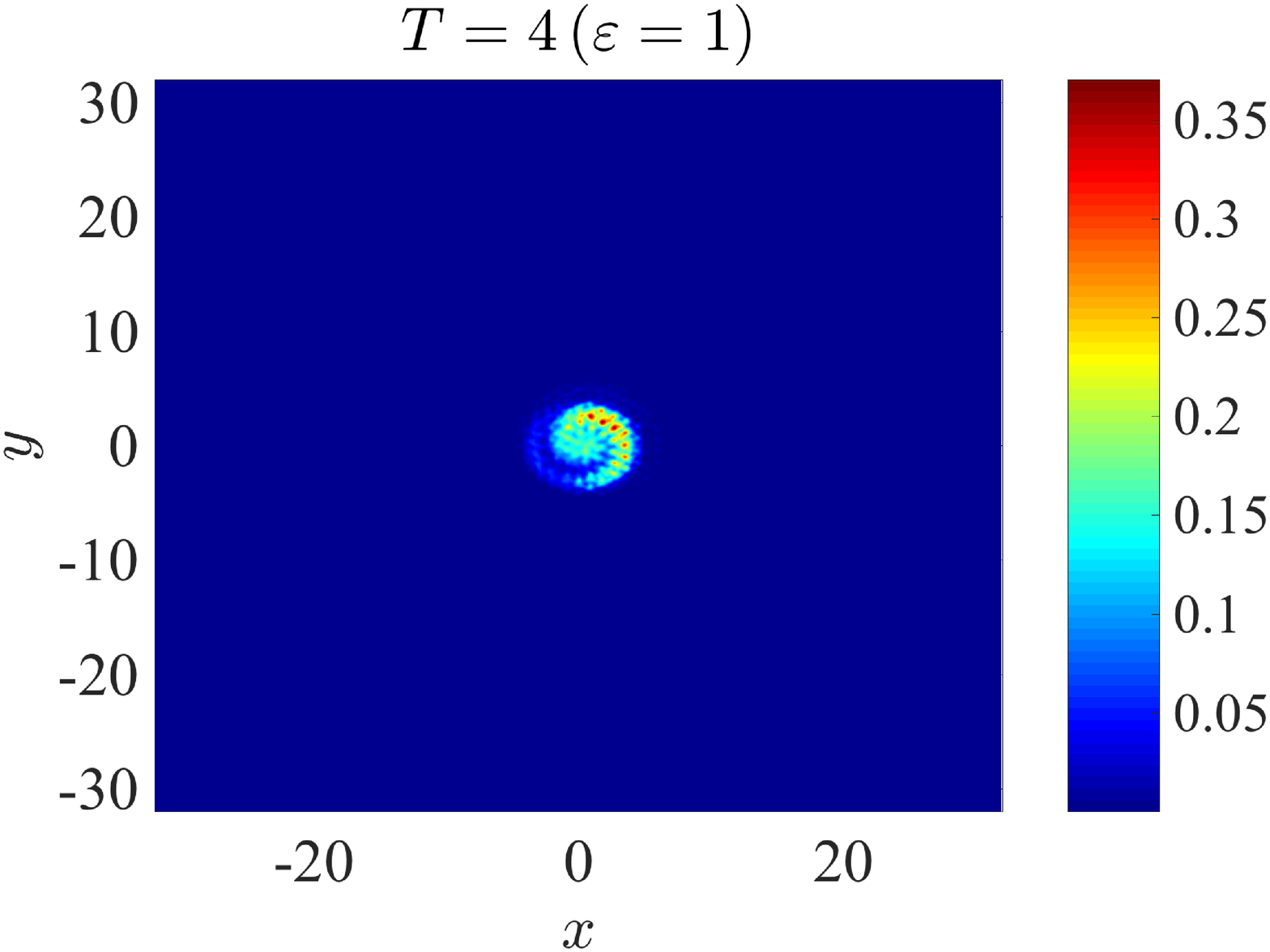}}
\end{minipage}
\begin{minipage}{0.31\textwidth}
\centerline{\includegraphics[width=4.3cm,height=3.6cm]{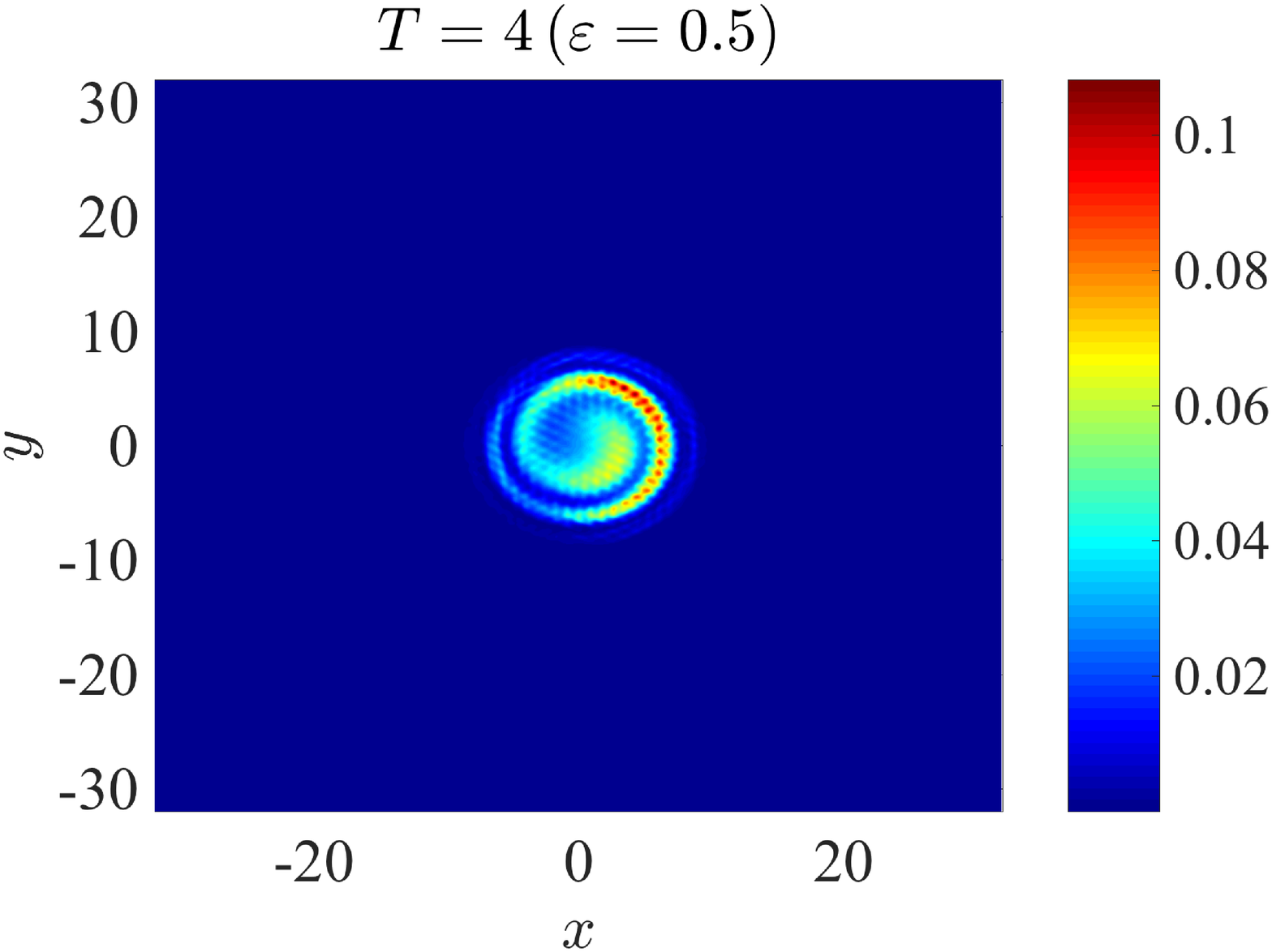}}
\end{minipage}
\begin{minipage}{0.31\textwidth}                                       \centerline{\includegraphics[width=4.3cm,height=3.6cm]{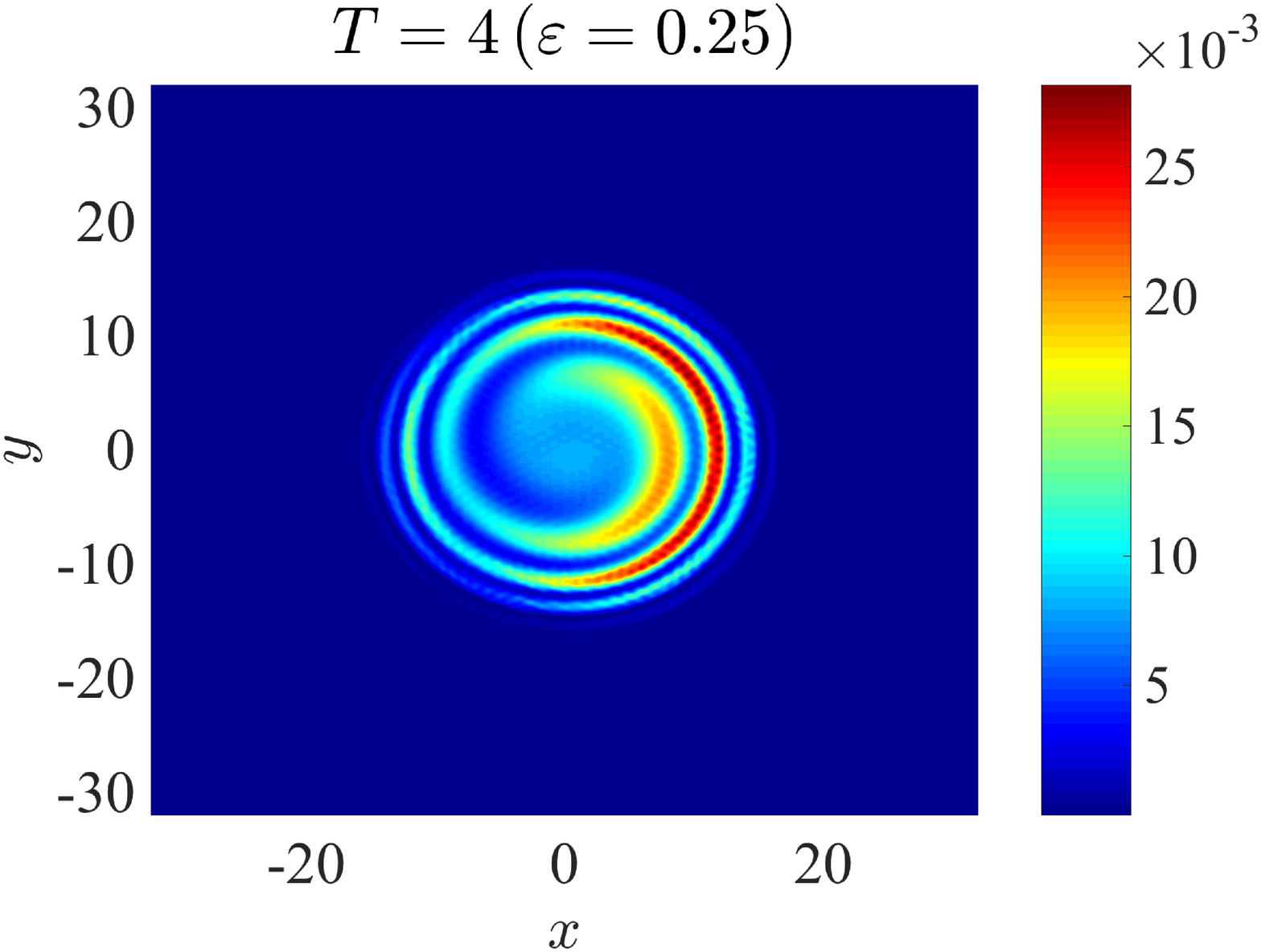}}
\end{minipage} \\
\begin{minipage}{0.31\textwidth}
\centerline{\includegraphics[width=4.3cm,height=3.6cm]{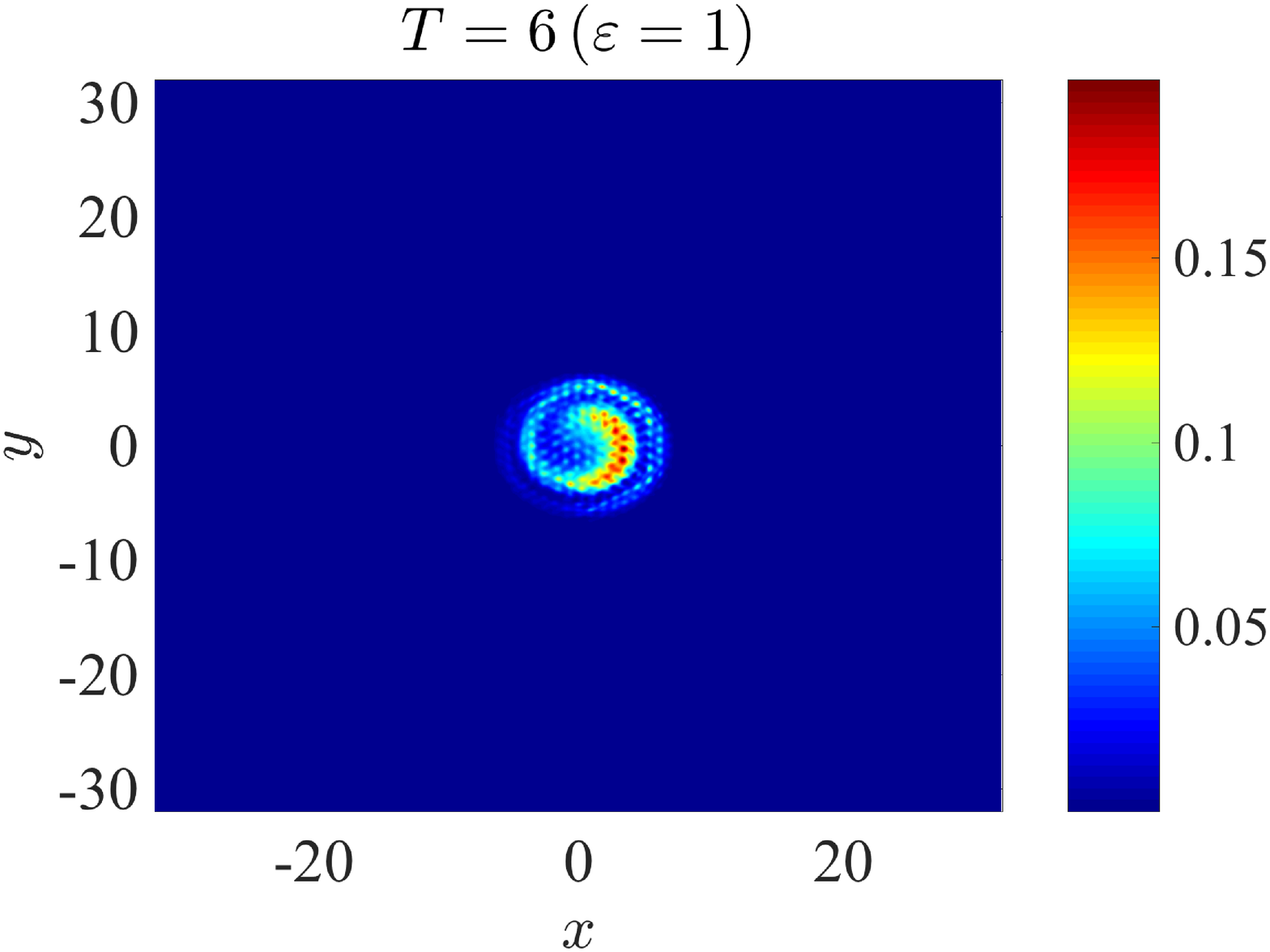}}
\end{minipage}
\begin{minipage}{0.31\textwidth}
\centerline{\includegraphics[width=4.3cm,height=3.6cm]{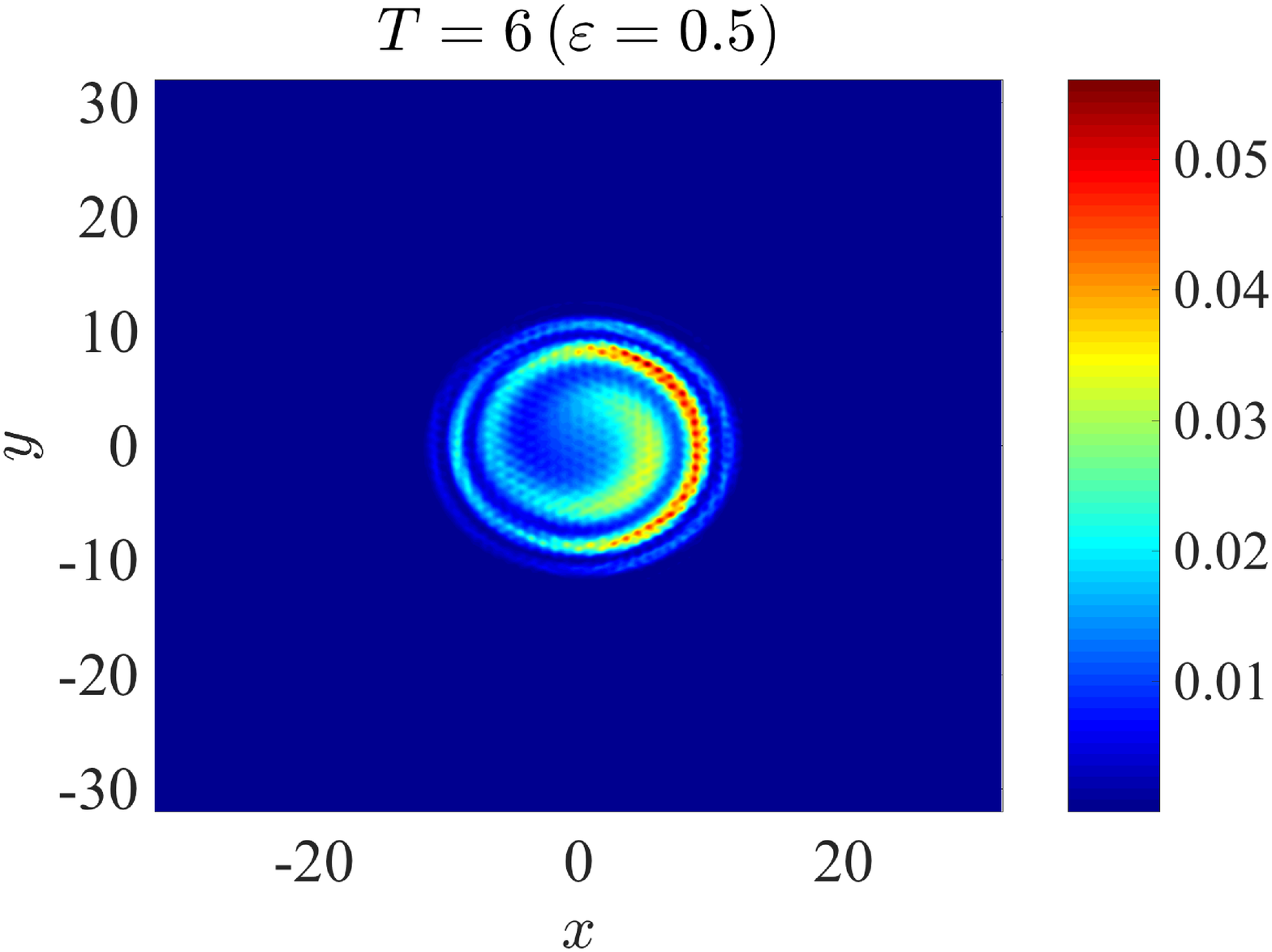}}
\end{minipage}
\begin{minipage}{0.31\textwidth}
\centerline{\includegraphics[width=4.3cm,height=3.6cm]{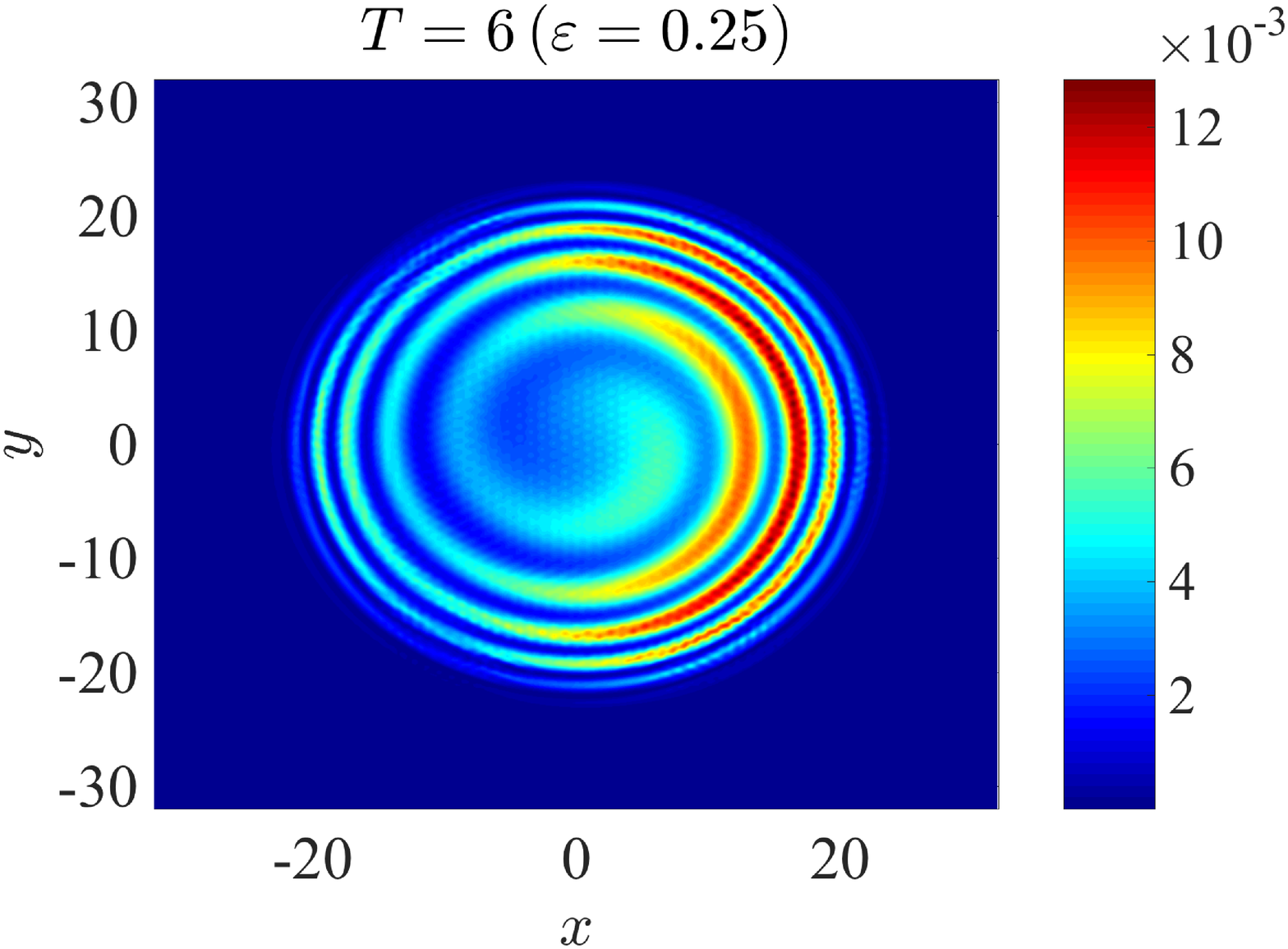}}
\end{minipage}  \\
\begin{minipage}{0.31\textwidth}
\centerline{\includegraphics[width=4.3cm,height=3.6cm]{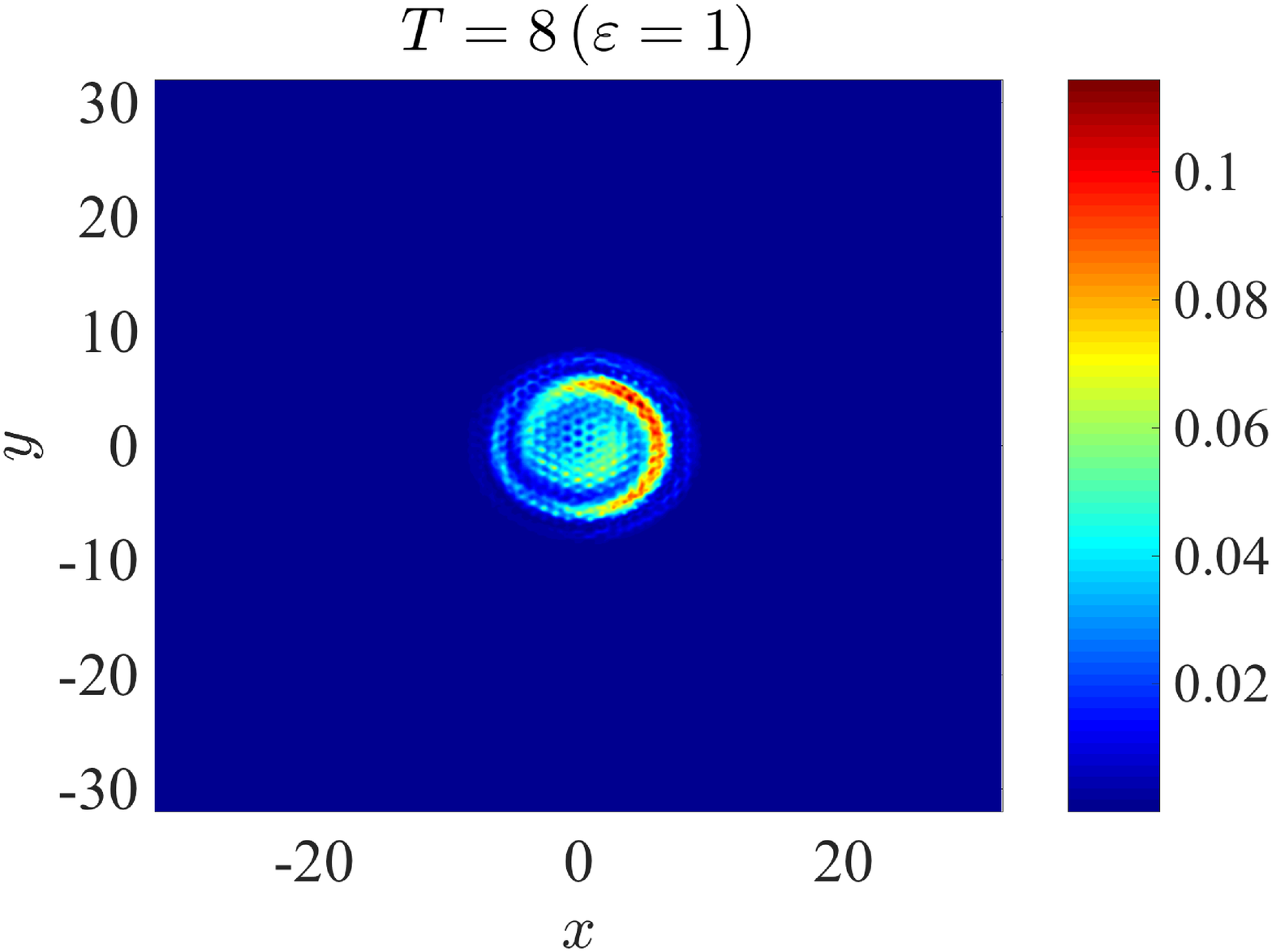}}
\end{minipage}
\begin{minipage}{0.31\textwidth}
\centerline{\includegraphics[width=4.3cm,height=3.6cm]{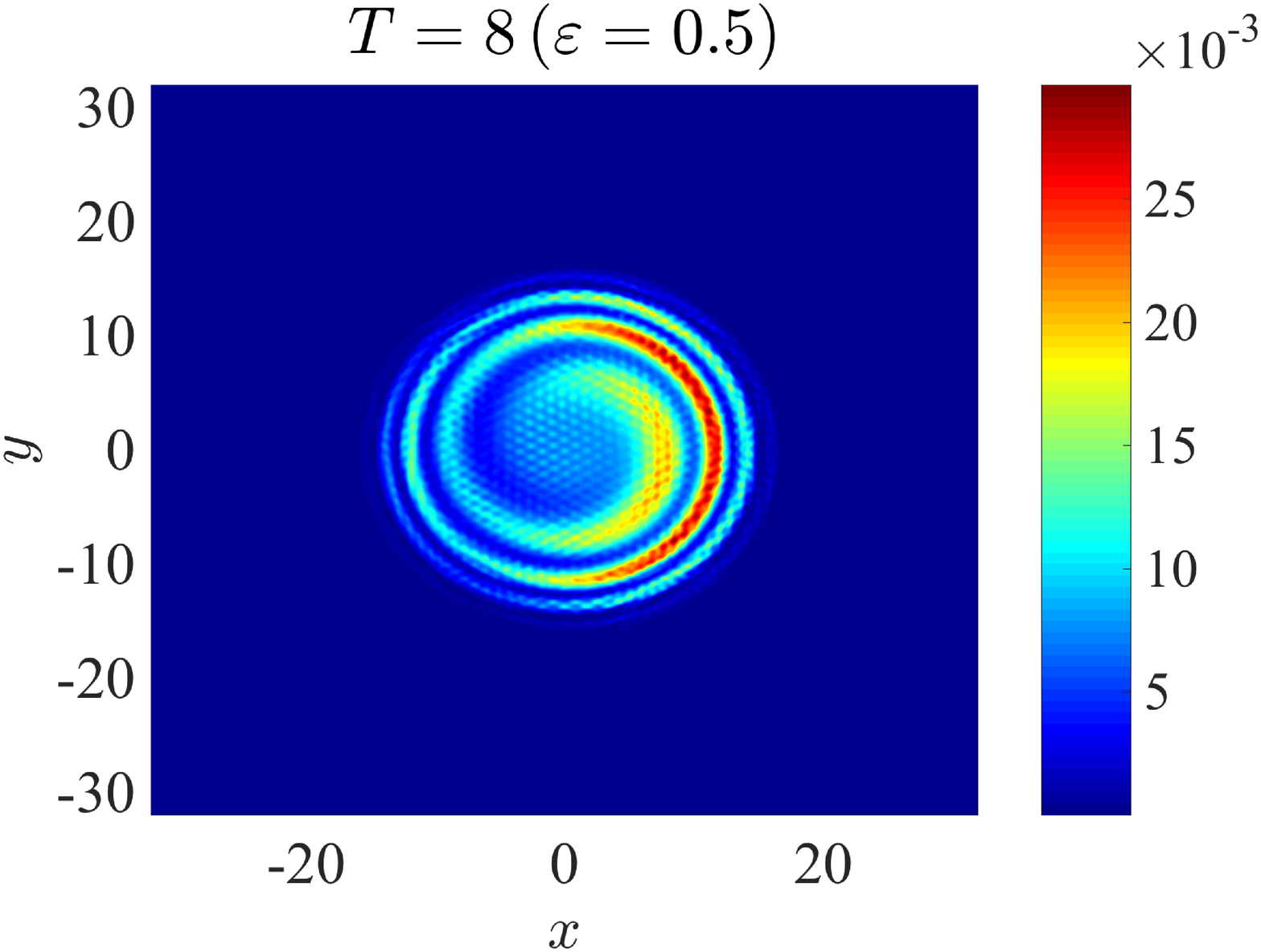}}
\end{minipage}
\begin{minipage}{0.31\textwidth}
\centerline{\includegraphics[width=4.3cm,height=3.6cm]{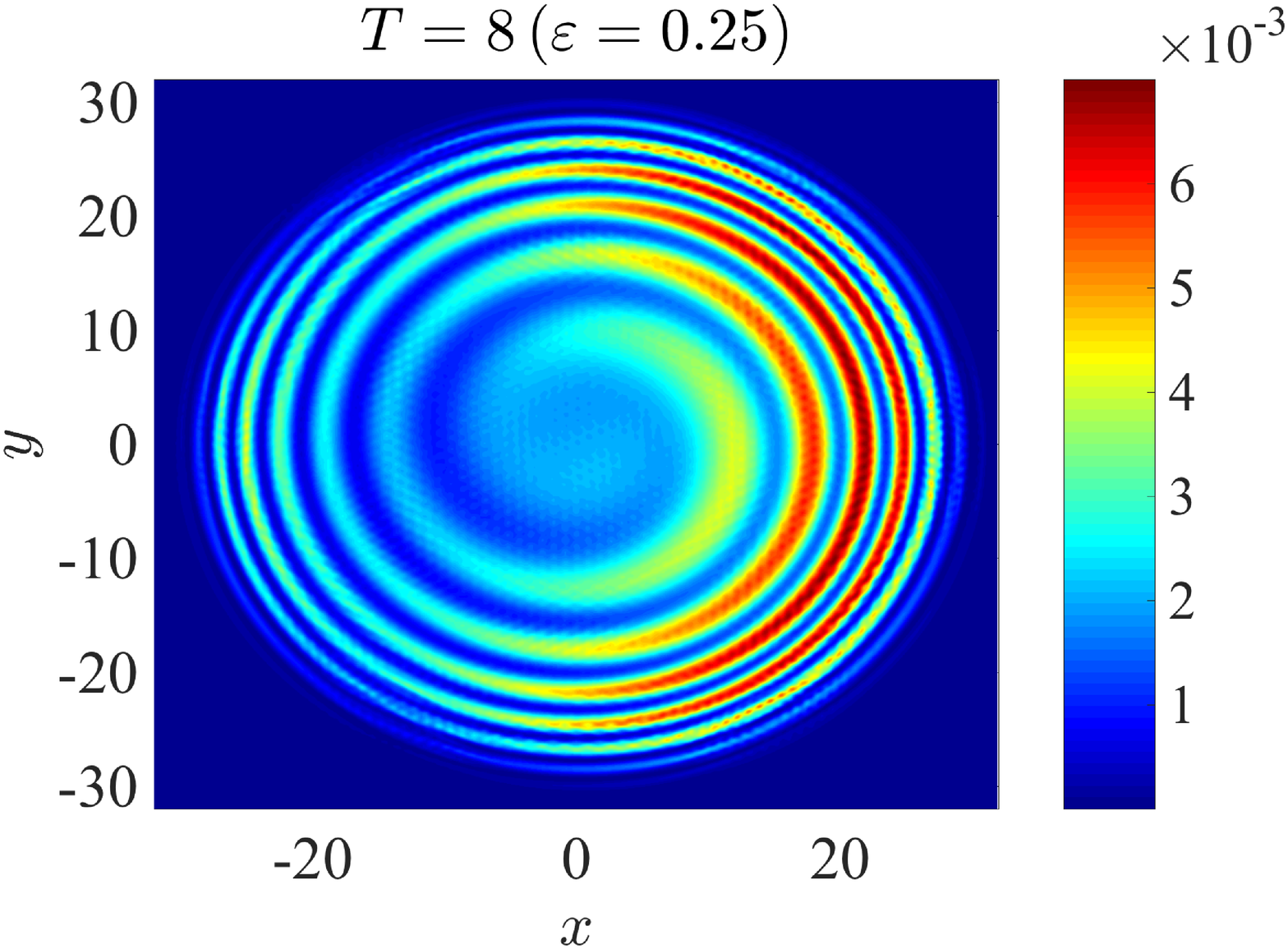}}
\end{minipage}
\caption{Dynamics of the total density $\rho$ of the Dirac equation \eqref{eq:Dirac_21} in 2D
with the honeycomb lattice potential \eqref{honpot} at different time $T$ for
different $\varepsilon$: $\varepsilon=1$ (left column),
$\varepsilon=0.5$ (middle column), and $\varepsilon=0.25$ (right column)}
\label{fig:2DLPe}
\end{figure}

{\textbf{Example 4.1}} (Honeycomb lattice potential).
We take $A_1(t,\textbf{x})=A_2(t,\textbf{x})\equiv0$ and a honeycomb lattice potential
\begin{equation}\label{honpot}
V(t,\textbf{x})=
\cos\left(\frac{4\pi}{\sqrt{3}}\textbf{e}_1\cdot\textbf{x}\right)+
\cos\left(\frac{4\pi}{\sqrt{3}}\textbf{e}_2\cdot\textbf{x}\right)+
\cos\left(\frac{4\pi}{\sqrt{3}}\textbf{e}_3\cdot\textbf{x}\right),
\end{equation}
with
\begin{equation}\label{honvec}
\textbf{e}_1=(-1,0)^T,\quad \textbf{e}_2=(1/2,\sqrt{3}/2)^T,\quad
\textbf{e}_3=(1/2,-\sqrt{3}/2)^T.
\end{equation}
The initial data is chosen as
\begin{equation}\label{honini}
\phi_1(0,\textbf{x})=e^{-\frac{x^2+y^2}{2}},\quad
\phi_2(0,\textbf{x})=e^{-\frac{(x-1)^2+y^2}{2}},\quad
\textbf{x}=(x,y)^T\in \Omega.
\end{equation}
The problem is solved numerically on $\Omega = (-32,32)^2$ by the
semi-implicit 4cFD method with the mesh size $h=1/16$ and time step $\tau=0.01$. Fig. \ref{fig:2DLPe} plots the total density $\rho(t,\textbf{x})$ of the Dirac equation
\eqref{eq:Dirac_21} at different time $T$ for different
$\varepsilon$.

From each column in Fig. \ref{fig:2DLPe}, we could observe that the dynamics of the Dirac equation \eqref{eq:Dirac_21} depends heavily on the parameter $\varepsilon$. For the classical regime, i.e., $\varepsilon= 1$, the density fluctuates in a random pattern as the stated Zitterbewegung for the relativistic dynamics (cf. left column in Fig. \ref{fig:2DLPe}). When $\varepsilon$ becomes smaller, the relativistic and mass effects are less and the density spreads over the lattice potential more smoothly and the maximal value becomes smaller (cf. each row in Fig. \ref{fig:2DLPe}). In addition, when $\varepsilon$ becomes half, the wave speed is double (cf. each row in Fig. \ref{fig:2DLPe} for different $\varepsilon$), which again confirms the solution propagates wave with wave speed at $O(1/\varepsilon)$. As a result, the 4cFD method could capture the dynamics of the Dirac equation accurately.

\bigskip

\begin{figure}[h!]
\begin{minipage}{0.31\textwidth}
\centerline{\includegraphics[width=4.3cm,height=3.6cm]{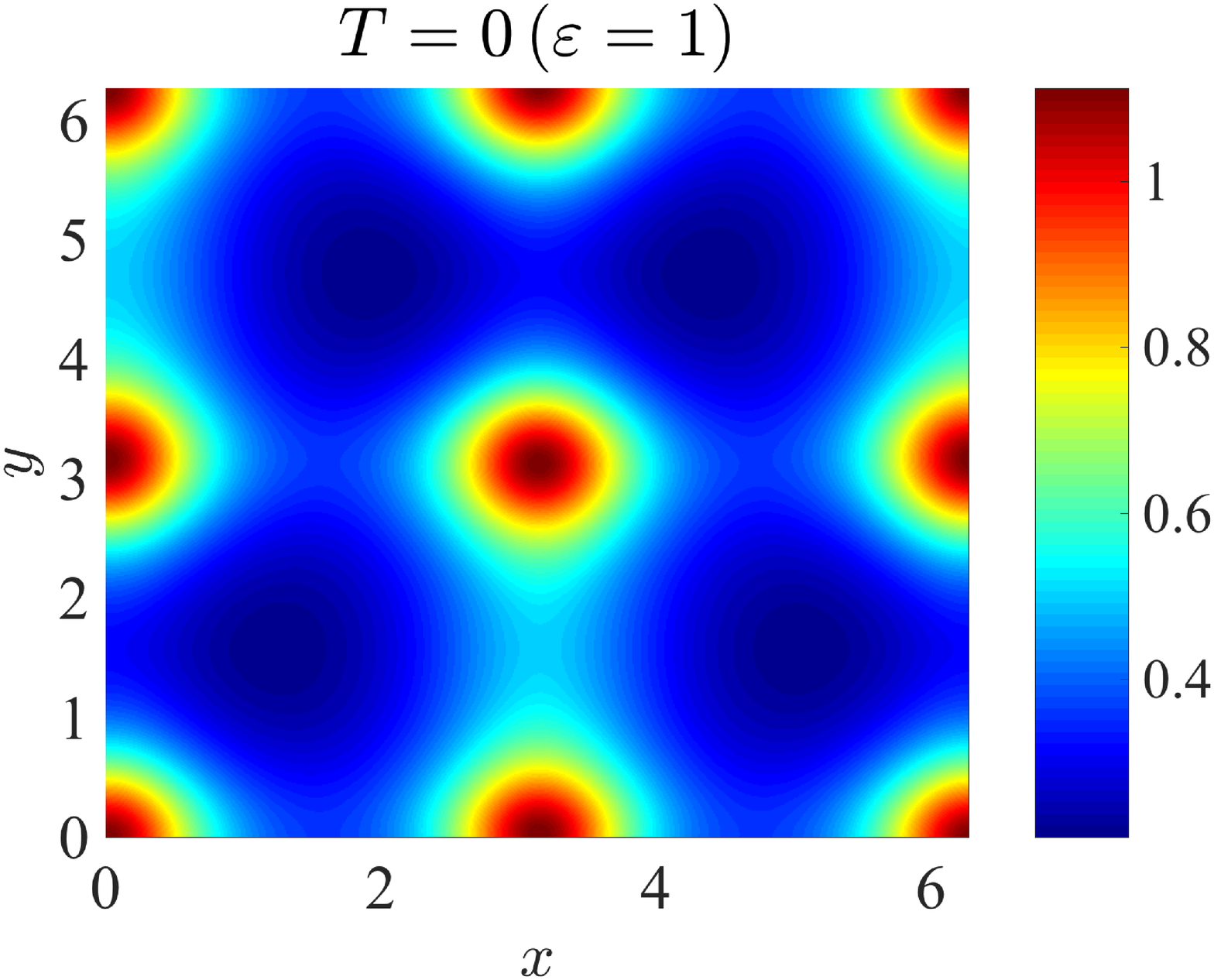}}
\end{minipage}
\begin{minipage}{0.31\textwidth}
\centerline{\includegraphics[width=4.3cm,height=3.6cm]{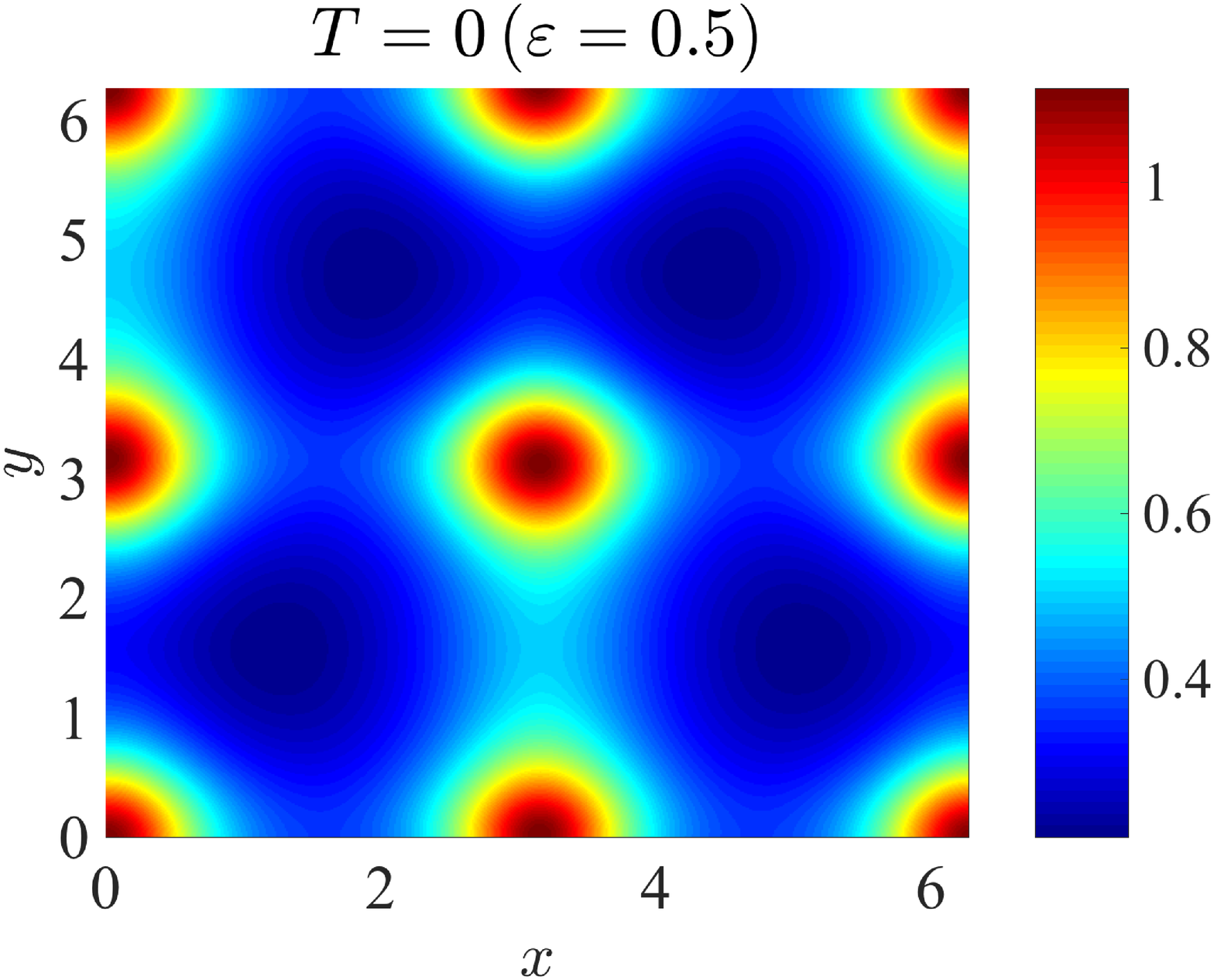}}
\end{minipage}
\begin{minipage}{0.31\textwidth}
\centerline{\includegraphics[width=4.3cm,height=3.6cm]{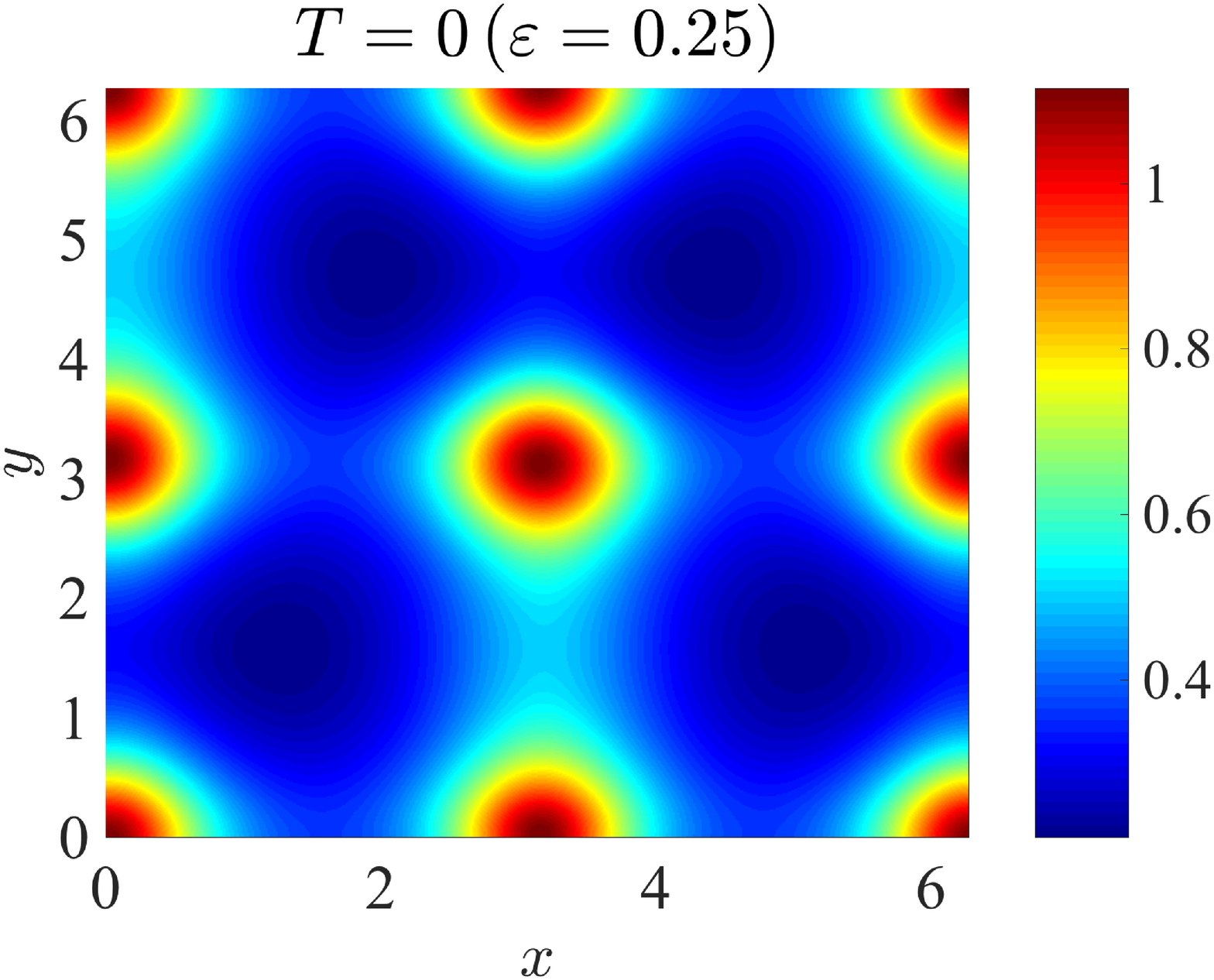}}
\end{minipage} \\
\begin{minipage}{0.31\textwidth}
\centerline{\includegraphics[width=4.3cm,height=3.6cm]{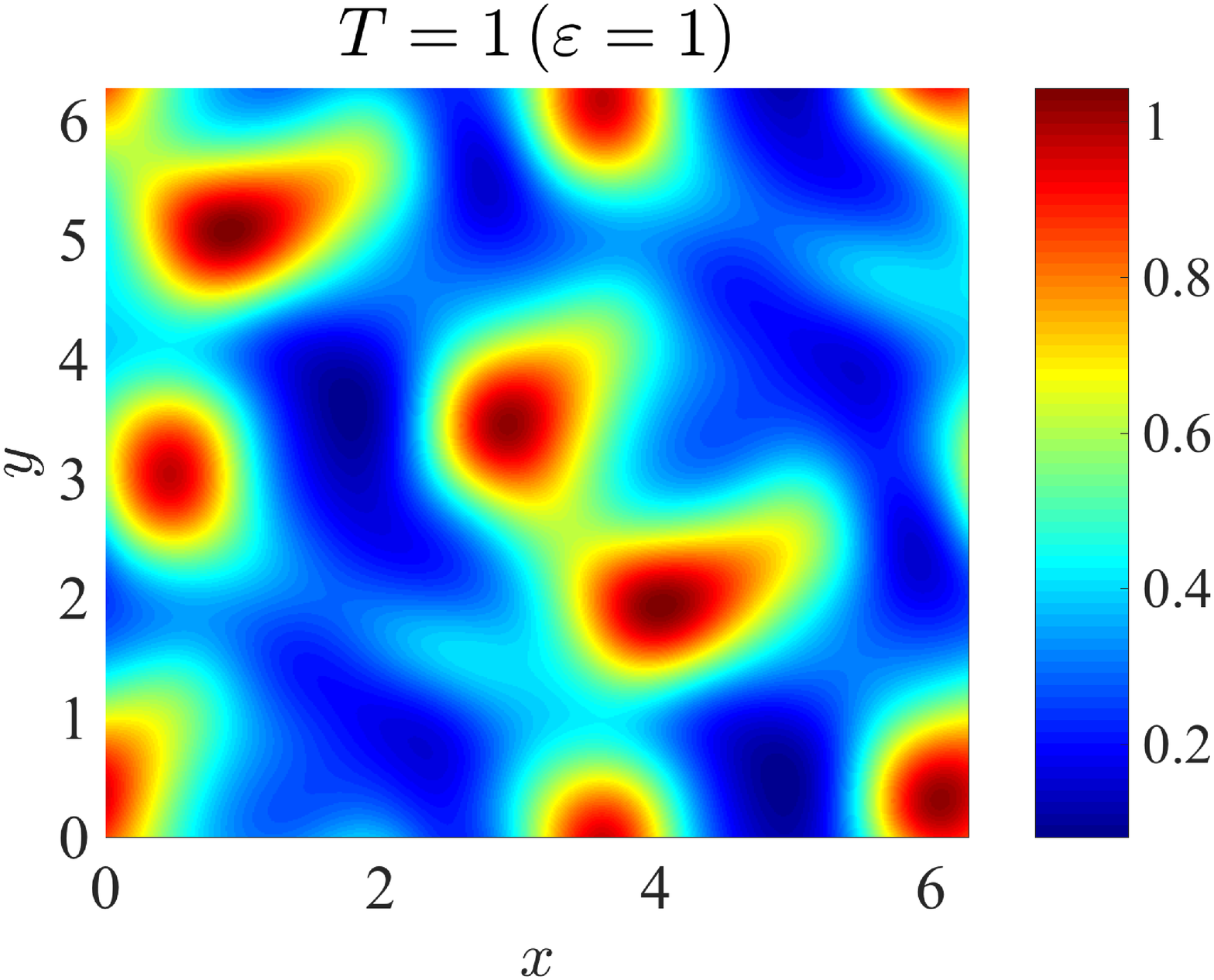}}
\end{minipage}
\begin{minipage}{0.31\textwidth}
\centerline{\includegraphics[width=4.3cm,height=3.6cm]{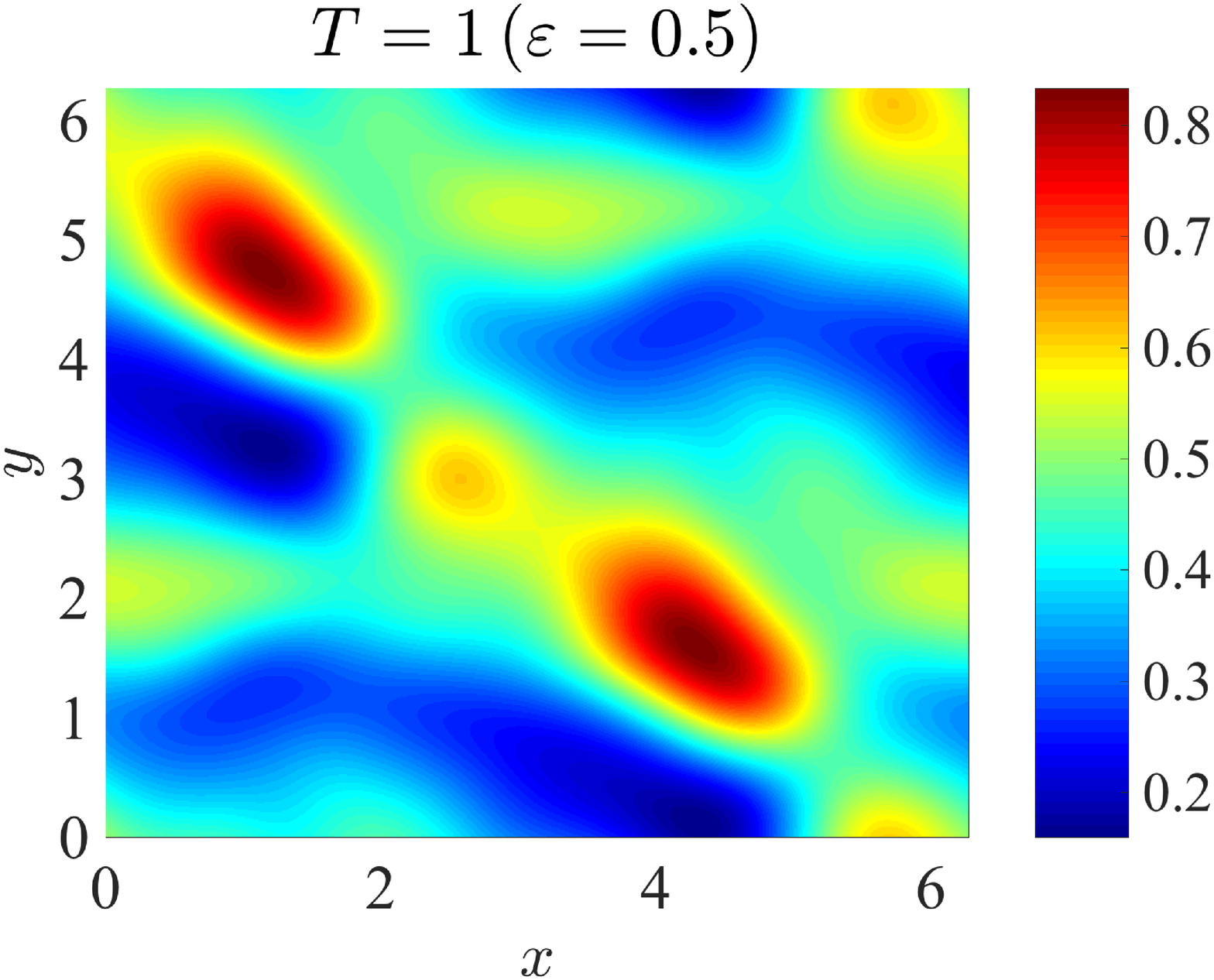}}
\end{minipage}
\begin{minipage}{0.31\textwidth}
\centerline{\includegraphics[width=4.3cm,height=3.6cm]{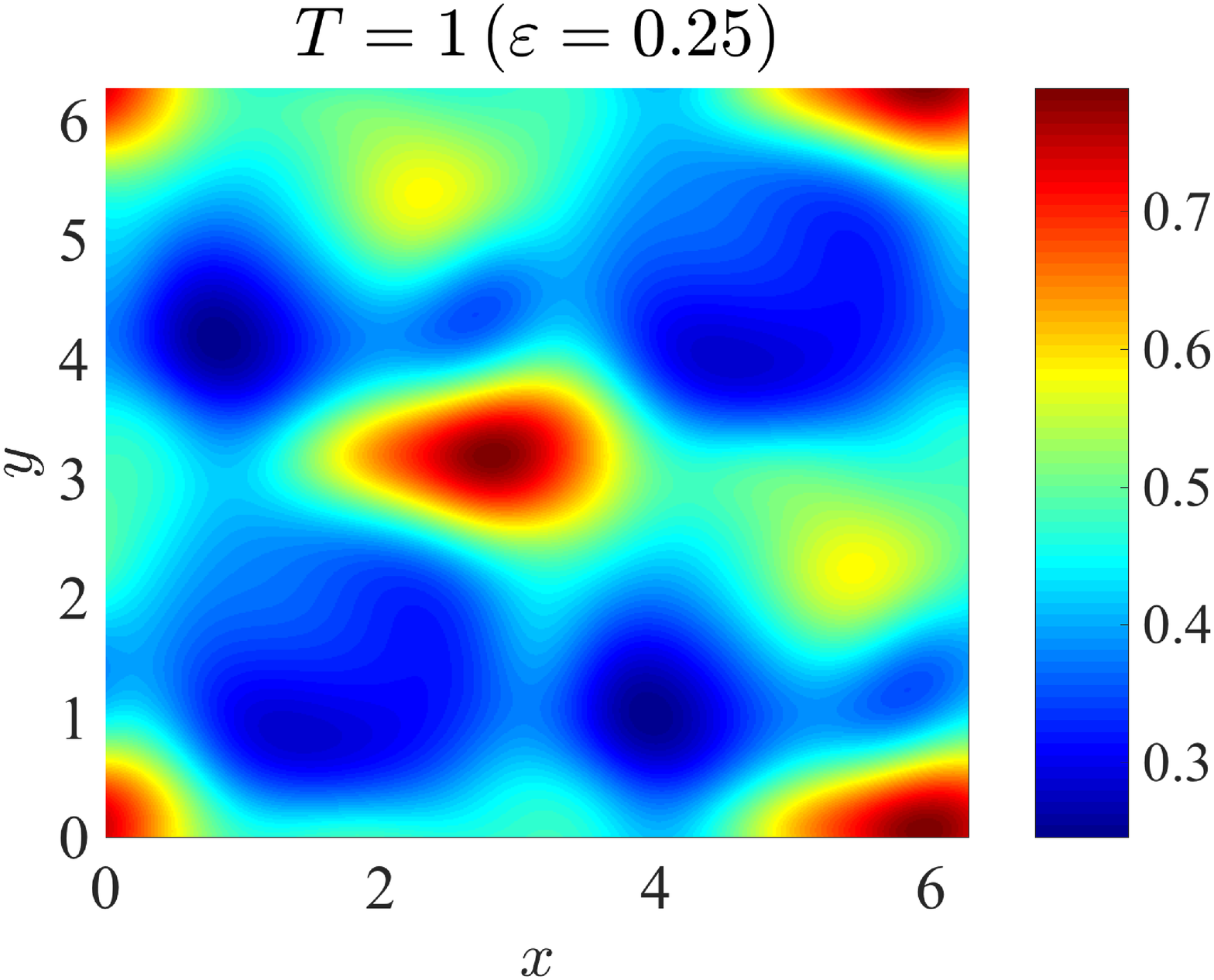}}
\end{minipage}\\
\begin{minipage}{0.31\textwidth}
\centerline{\includegraphics[width=4.3cm,height=3.6cm]{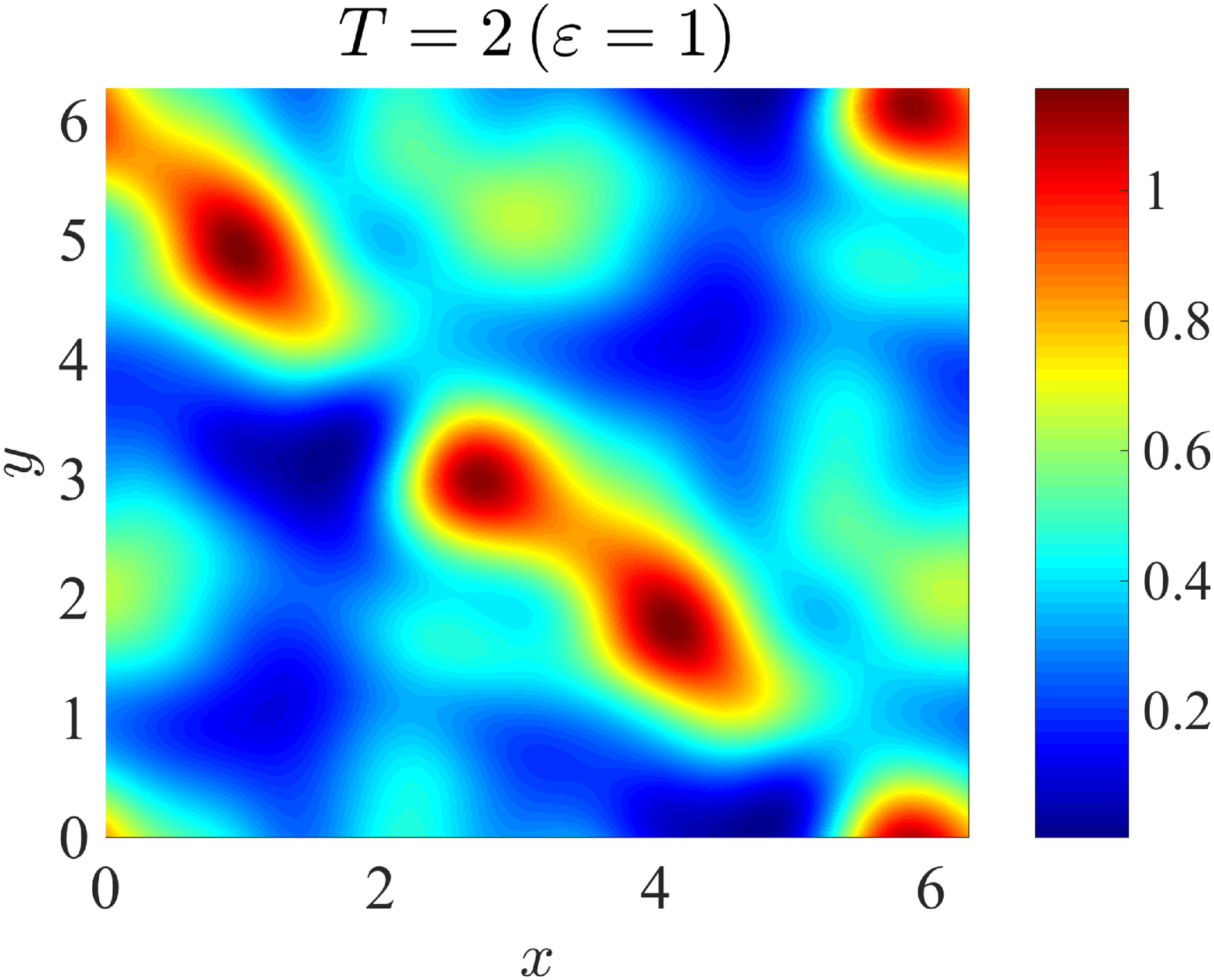}}
\end{minipage}
\begin{minipage}{0.31\textwidth}
\centerline{\includegraphics[width=4.3cm,height=3.6cm]{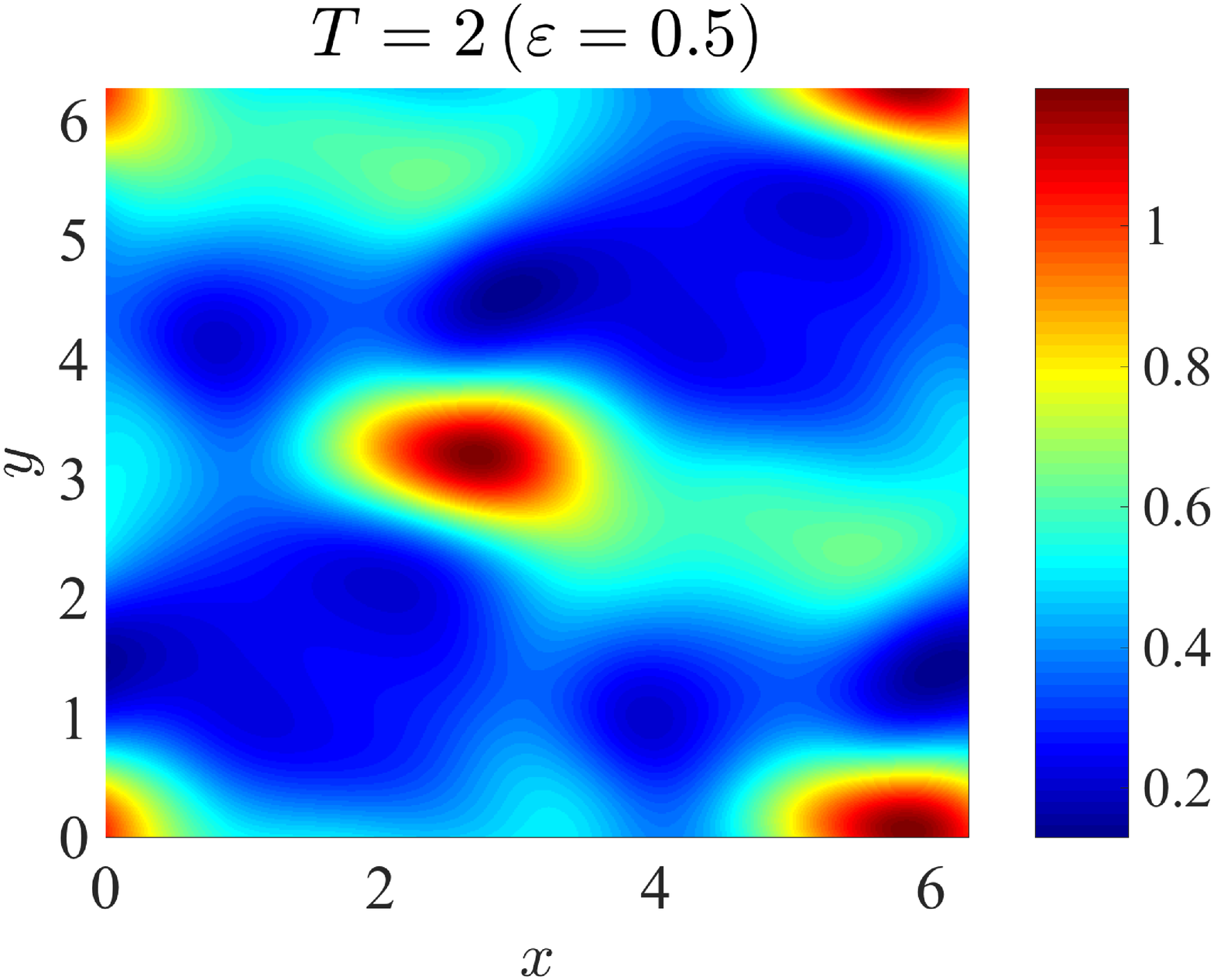}}
\end{minipage}
\begin{minipage}{0.31\textwidth}
\centerline{\includegraphics[width=4.3cm,height=3.6cm]{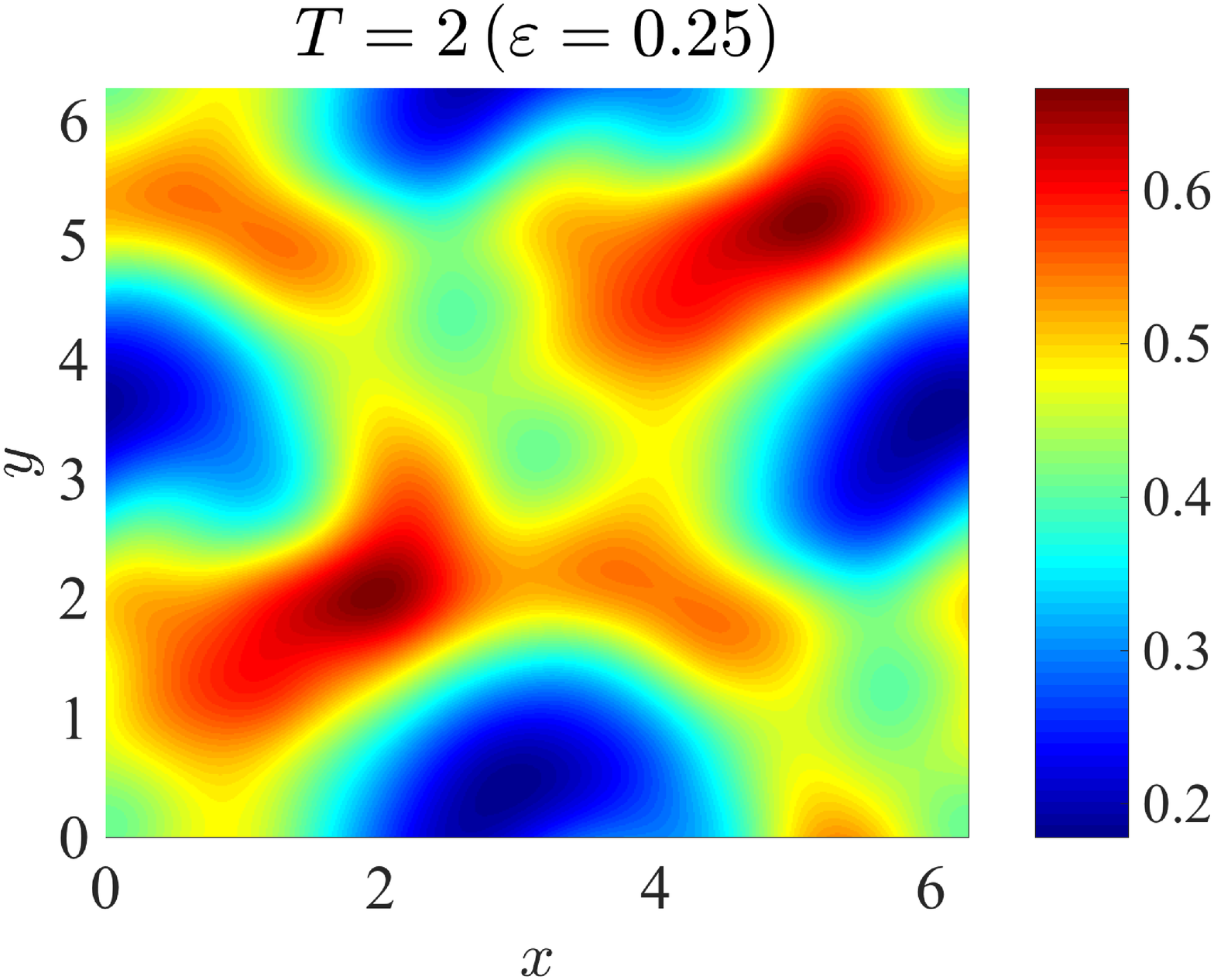}}
\end{minipage} \\
\begin{minipage}{0.31\textwidth}
\centerline{\includegraphics[width=4.3cm,height=3.6cm]{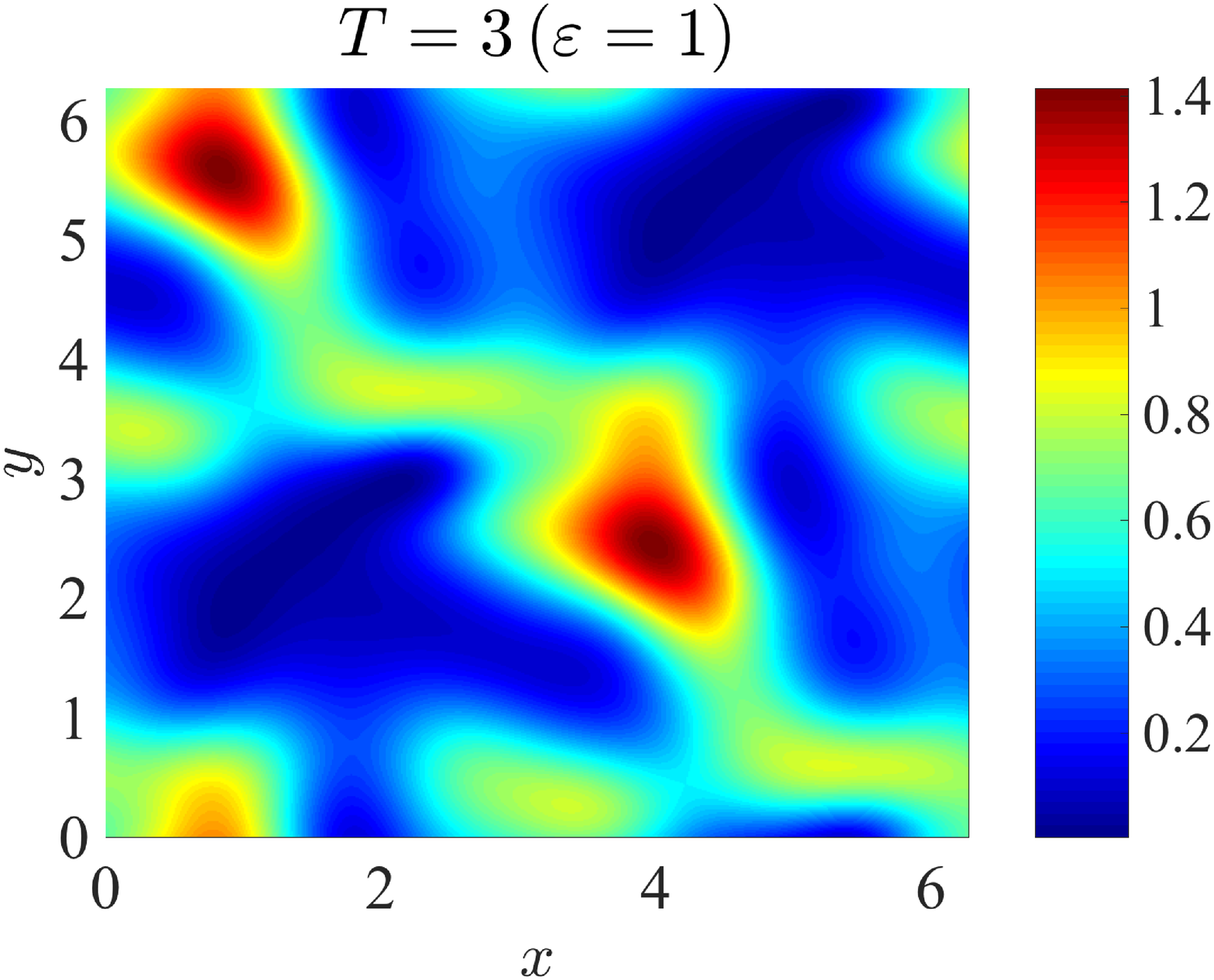}}
\end{minipage}
\begin{minipage}{0.31\textwidth}
\centerline{\includegraphics[width=4.3cm,height=3.6cm]{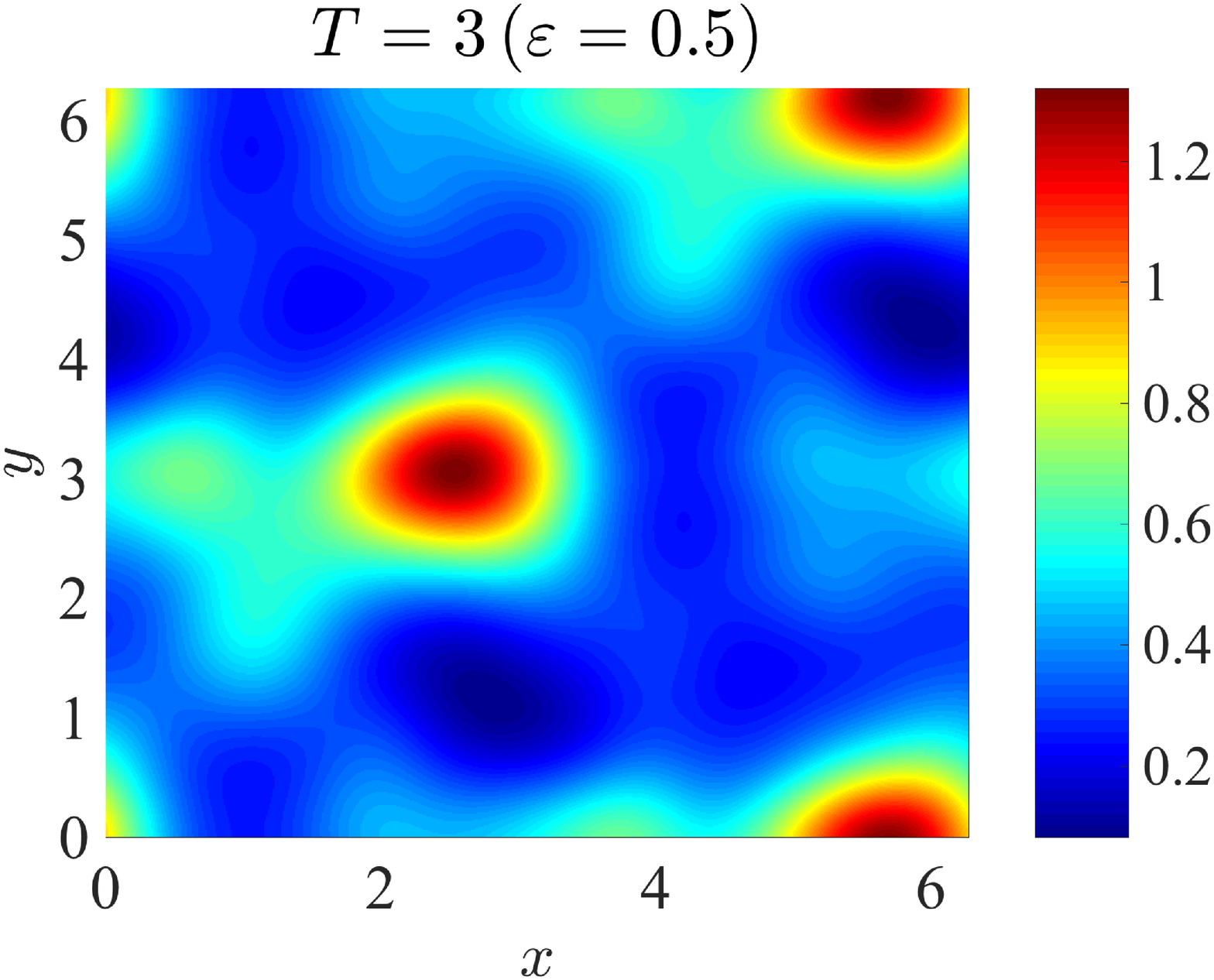}}
\end{minipage}
\begin{minipage}{0.31\textwidth}
\centerline{\includegraphics[width=4.3cm,height=3.6cm]{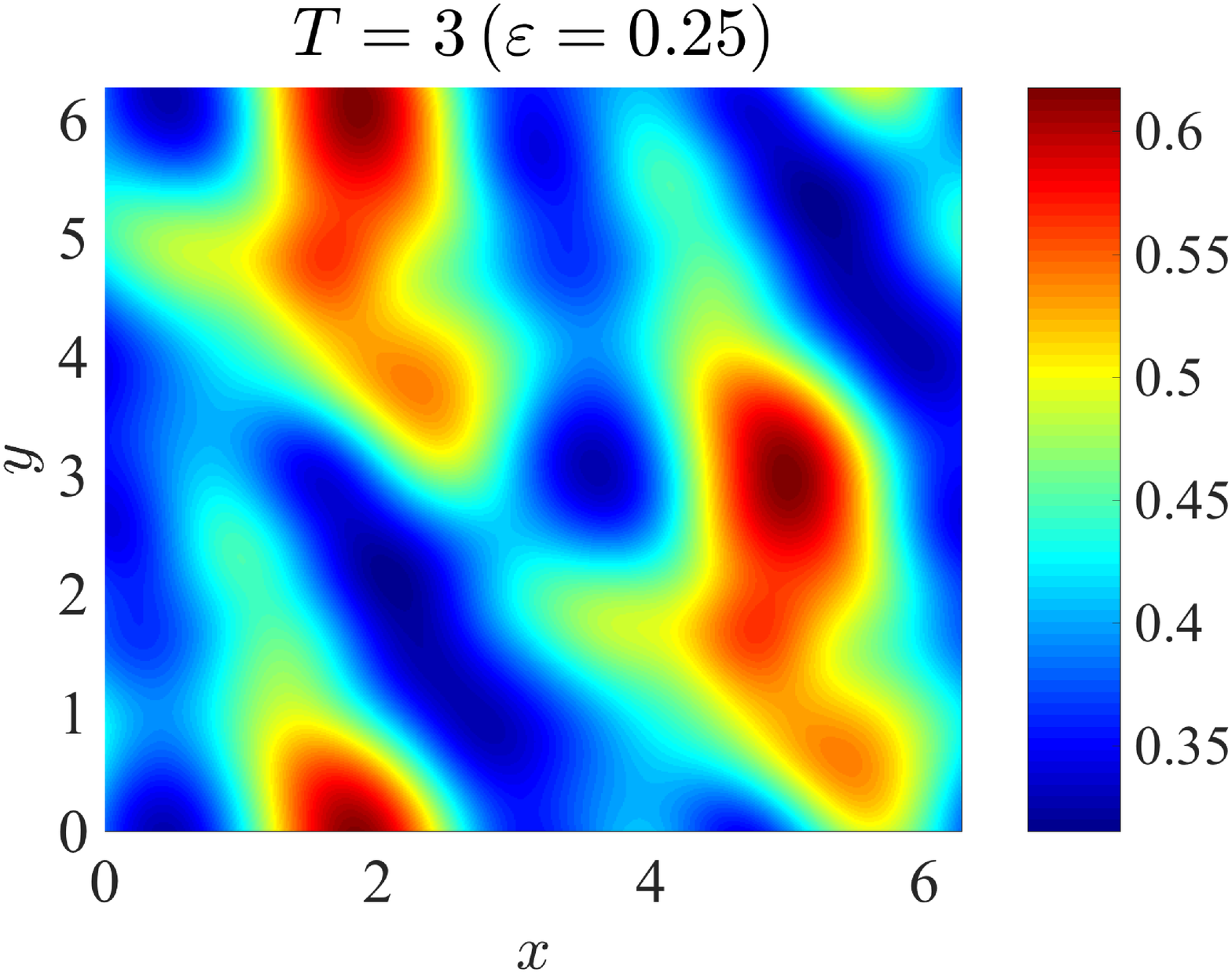}}
\end{minipage}  \\
\begin{minipage}{0.31\textwidth}
\centerline{\includegraphics[width=4.3cm,height=3.6cm]{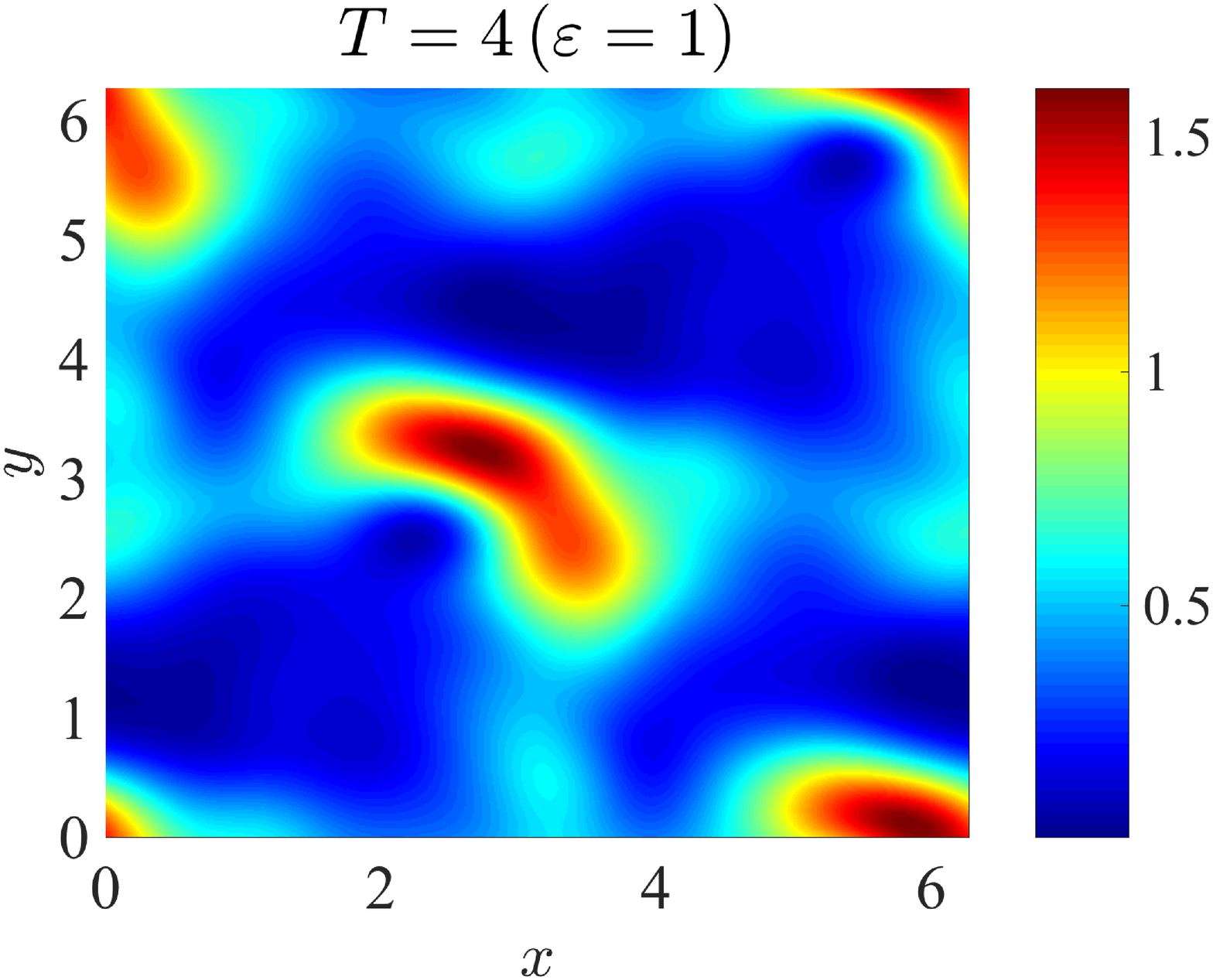}}
\end{minipage}
\begin{minipage}{0.31\textwidth}
\centerline{\includegraphics[width=4.3cm,height=3.6cm]{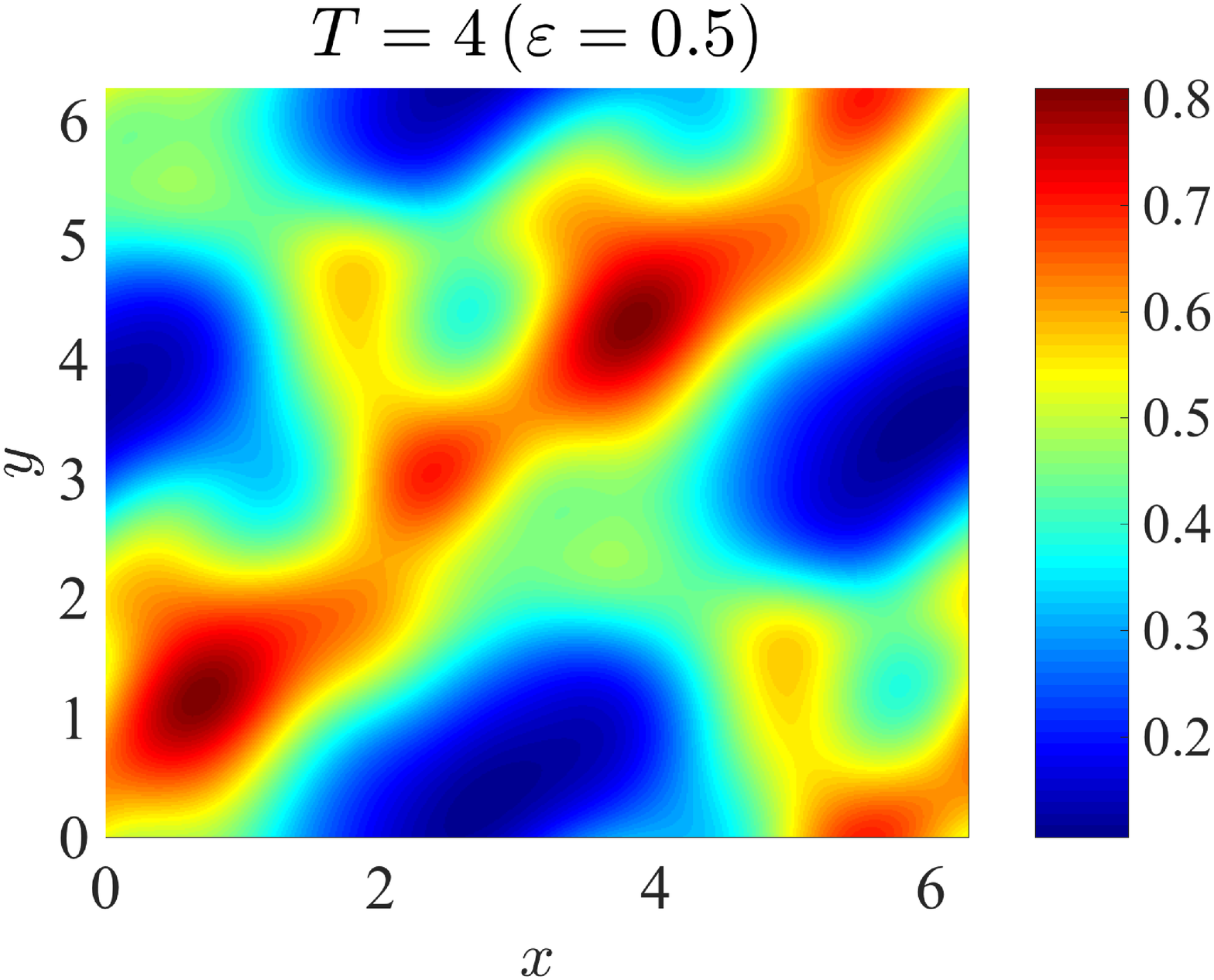}}
\end{minipage}
\begin{minipage}{0.31\textwidth}
\centerline{\includegraphics[width=4.3cm,height=3.6cm]{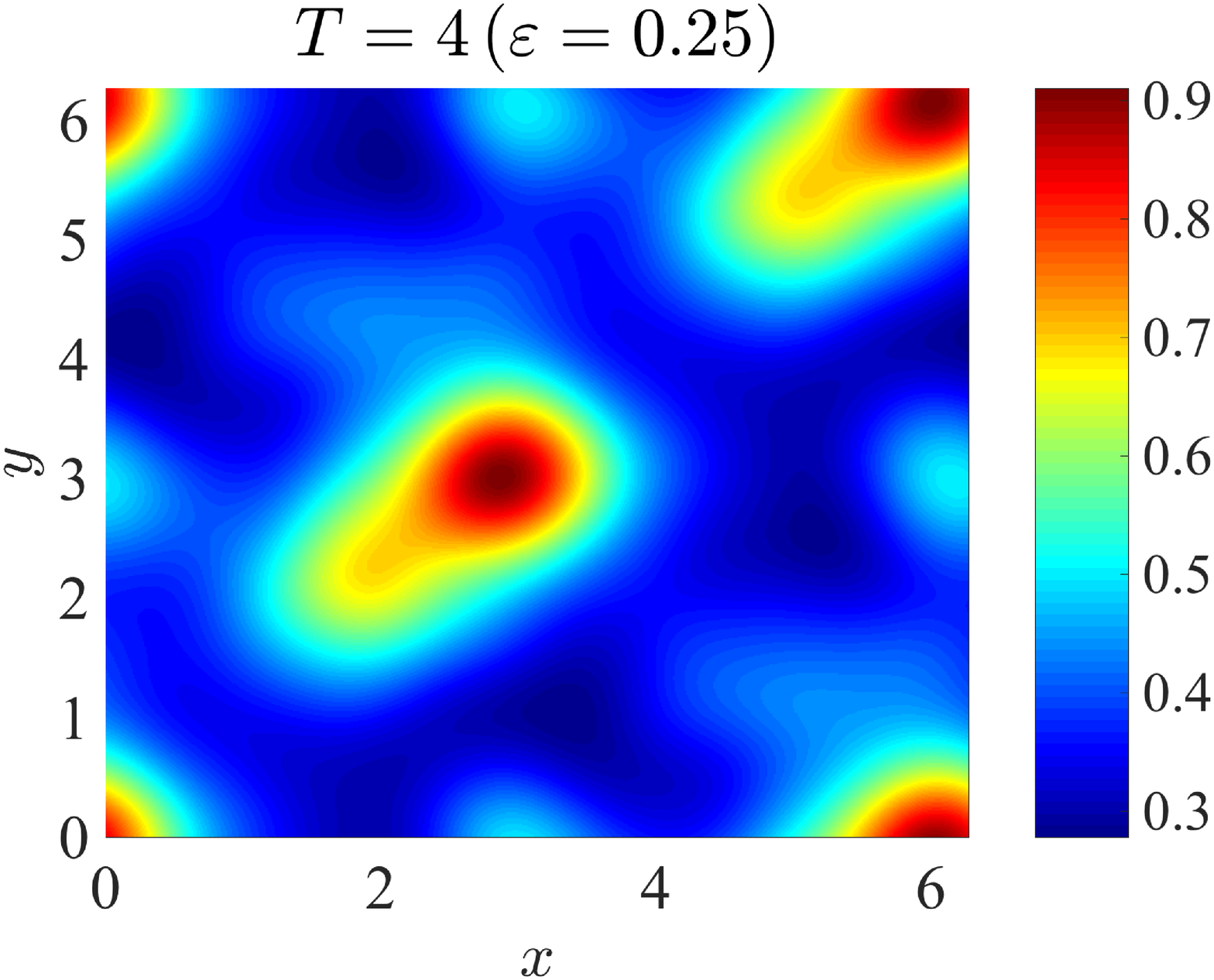}}
\end{minipage}
\caption{Dynamics of the total density $\rho$ of the Dirac equation \eqref{eq:Dirac_21} in 2D
with the periodic electromagnetic potentials \eqref{empot1} at different time $T$ for
different $\varepsilon$: $\varepsilon=1$ (left column),
$\varepsilon=0.5$ (middle column), and $\varepsilon=0.25$ (right column)}
\label{fig:2DGPe}
\end{figure}

{\textbf{Example 4.2}} (Periodic electromagnetic potentials). Here we take the following periodic electromagnetic potentials
\begin{equation}\label{empot1}
V(t, \textbf{x}) =   \frac{1}{1+\sin^2(x)+\cos^2(y)}, \quad
A_1(t,\textbf{x})= \sin(2x)\sin(2y),\quad A_2(t,\textbf{x}) = 2\sin(x)\cos(y),
\end{equation}
and the initial data as
\begin{equation}
\phi_1(0, \textbf{x}) = \frac{1}{1+\sin^2(x)+\sin^2(y)},\quad \phi_2(0, \textbf{x}) = \frac{1}{3+\cos(x)\sin(y)},\quad \textbf{x}=(x,y)^T\in\Omega.
\end{equation}
The problem is solved on $\Omega = (0, 2\pi)^2$ by the semi-implicit 4cFD method with the mesh size $h=\pi/256$ and time step $\tau=0.01$. Fig. \ref{fig:2DGPe} depicts the total density $\rho(t,\textbf{x})$ of the Dirac equation \eqref{eq:Dirac_21} at different time $T$ for different $\varepsilon$.

From each column in Fig. \ref{fig:2DGPe}, we find that the dynamics of the Dirac equation with periodic electromagnetic potentials also varies considerably for different $\varepsilon$. Again, the numerical results show that the 4cFD method could capture the dynamics very accurately.

\section{Conclusions}
\label{sec:5}
In this paper, the conservative and semi-implicit fourth-order compact finite difference (4cFD) methods were applied to numerically solve the Dirac equation in the massless and nonrelativistic regime. The mass and energy conservation, stability conditions and error bounds of the 4cFD methods for the Dirac equation were rigorously established. The error bounds depend explicitly on the mesh size $h$ and time step $\tau$ as well as the small parameter $\varepsilon \in (0, 1]$ and indicate that the $\varepsilon$-scalability of the 4cFD  methods should be taken as $h=O(\varepsilon^{1/4})$ and $\tau = O(\varepsilon^{3/2})$, which have higher spatial convergence and better spatial resolution than the classical finite difference methods. Furthermore, the physical observations including the the total density and current density also have the same resolution in this regime. Numerical results were presented to confirm the error bounds and the $\varepsilon$-scalability of the 4cFD methods. In addition, the semi-implicit scheme was used to study the dynamics of the Dirac equation with different electromagnetic potentials in 2D and some interesting phenomena were observed. \\

\section*{Acknowledgements}

This work was supported by the Ministry of Education of Singapore Grant No. R-146-000-290-114 (Y. Feng), and the National Natural Science Foundation of China Grant No. U1930402 (Y. Ma).


\end{document}